\documentstyle[11pt,amssymb,xypic]{article}

\newcommand{\sq}{{$\enspace\square$}}
\newcommand{\into}{\hookrightarrow}
\newcommand{\sto}{\rightsquigarrow}
\newcommand{\ad}{{\mbox{Ad}\,}}

\newcommand{\RR}{\mbox{\sf M}}
\newcommand{\g}{{\frak g}}

\newcommand{\too}{\;\longrightarrow\;}
\newcommand{\onto}{\;\twoheadrightarrow\;}
\newcommand{\st}{\divideontimes}
\newcommand{\alt}{{\bold\varepsilon}}
\newcommand{\mod}{\mbox{-mod}}

\newcommand{\La}{{\bold \Lambda}}
\newcommand{\eps}{{\mathbf{q}}}
\newcommand{\zht}{{\mathbf{t}}}
\newcommand{\epq}{q}
\newcommand{\CqepsT}{{{\C}(T)}}

\newcommand{\Waff}{{\vphantom{W}\smash{\buildrel_{\,\,{_{_\bullet}}}\over W}}}
\newcommand{\MM}{{\mathfrak{M}}}
\newcommand{\MMM}{{\mathbf{M}}}
\newcommand{\m}{{\mathbf{m}}}

\newcommand{\dd}{{\mathsf{D}}_{_{^{\mathbf{q}}}}}
\newcommand{\ccirc}{{}_{^\circ}}
\newcommand{\TKM}{{\tilde{T}}}
\newcommand{\TKMaqm}{{\TKM}_{\alpha,\mathbf{t}^{-2}}}

\newcommand{\TKMgamma}{{\TKM}_{\gamma}}
\newcommand{\TKMgqm}{{\TKM}_{\gamma,\mathbf{t}^{-2}}}
\newcommand{\Lie}{{\sl {Lie}}^{\,}}
\newcommand{\rhocheck}{\check{\rho}}

\newcommand{\Raf}{{\vphantom{\Delta}\smash{\buildrel_{\,\,{_\bullet}}\over \Delta}}}

\newcommand{\Raffreal}{{
{\vphantom{\Delta}\smash{\buildrel_{\,\,{_\bullet}}\over \Delta}}_{_{^{\mathsf re}}}}}

\newcommand{\Raffrepos}{\displaystyle
{{\vphantom{\Delta}\smash{\buildrel_{\,\,{_\bullet}}\over \Delta}}_{_{{^{\mathsf
re}}}}^{^{_+}}}}
\newcommand{\dT}{{\vphantom{T}\smash{\buildrel_{\,\,
{_{_{_\bullet}}}}\over T}}}
\newcommand{\T}{{\mathsf{T}}}
\newcommand{\hT}{{\mathsf{\hat{T}}}}

\newcommand{\hh}{{\cal H}}
\newcommand{\os}{{\cal M}^{^{\mbox{{\tiny ss}}}}\!}
\newcommand{\A}{{\cal A}}

\newcommand{\bb}{{\cal B}}

\newcommand{\Y}{{\bold Y}}
\newcommand{\C}{{\Bbb C}}
\newcommand{\Z}{{\Bbb Z}}

\newcommand{\TT}{{\Bbb T}}
\newcommand{\W}{{\Bbb W}}
\newcommand{\oo}{{\cal M}}

\newcommand{\inv}{^{-1}}
\newcommand{\iso}{\,{\stackrel{\sim}{\longrightarrow}\,}}

\newcommand{\ab}{{\hskip 8mm}}

\newcommand{\ee}{{\bold e}}
\newcommand{\V}{{\bold V}}
\newcommand{\la}{\lambda}
\newcommand{\vi}{{\sf {(i)}}$\;\;$}
\newcommand{\vii}{{\sf {(ii)}}$\;\;$}
\newcommand{\viii}{{\sf {(iii)}}$\;\;$}
\newcommand{\dis}{\displaystyle}

\newcommand{\HH}{{\mathsf H}}

\newcommand{\hdd}{{\vphantom{\footnotesize{L}}
\smash{\buildrel_{\,_{_{_{\bullet\bullet}}}}\over {\mathsf H}}}}
\newcommand{\hd}{{\vphantom{\footnotesize{L}}\smash{\buildrel_{\,{_\bullet}}\over
{\mathsf H}}}}

\newcommand{\ehe}{{\ee\!\hdd\!\ee}}

\newcommand{\hddv}{\hdd_{\!v}}

\newcommand{\ehev}{\ee_v\hddv \ee_v}
\newcommand{\ehez}{\ee_{_0}\hdd_{_0}\ee_{_0}}

\newcommand{\Hom}{{\mathsl{Hom}}}
\newcommand{\ind}{{\mathsl{Ind}\,}}
\newcommand{\Res}{{\mathsl{Res}\,}}

\newcommand{\hw}{{\cal H}^W}
\newcommand{\X}{{\bold X}}

\newcommand{\Higgs}{{\mathsf{Higgs}}}
\newcommand{\sset}{\subset}

\newcommand{\E}{{\mathcal{E}}}

\begin{document}

\setlength{\parindent}{0pt}
\setlength{\parskip}{3pt plus 5pt minus 1pt}

\centerline{\huge {\bf Representations of quantum tori}}
\vskip 3mm
\centerline{\huge {\bf  and double-affine Hecke algebras.}}

\vskip 10mm
\centerline{\large {\sc Vladimir Baranovsky, Sam Evens, and Victor Ginzburg}}
\vskip 2pt

\medskip
\begin{abstract}{\footnotesize
We study a BGG-type category of infinite dimensional
representations of $\hh[W]$, a semi-direct product of
the quantum torus with parameter $\eps$,
 built on the root lattice of a semisimple group $G$,
and the Weyl group of $G$. Irreducible objects of our category turn
out to be parameterized by semistable
$G$-bundles on the elliptic curve
$\C^*\!/\eps^\Z$. 
In the second part of the paper we construct a 
family of algebras depending on a parameter $v$ that specializes
to $\hh[W]$ at $v=0$, and specializes
to the double-affine Hecke algebra $\hdd$, introduced by
 Cherednik, at $v=1$. We propose a
Deligne-Langlands-Lusztig type
conjecture relating irreducible  $\hdd$-modules
to  Higgs
$G$-bundles on the elliptic curve.
 The conjecture may be seen as a natural `$v$-deformation' of
the classification of  simple $\hh[W]$-modules obtained in the
first part of the paper.
Also, an `operator realization'  of the double-affine Hecke algebra,
 as well as of its {\it Spherical subalgebra},
in terms of certain
`zero-residue' conditions
is given.}
\end{abstract}

\centerline{\sf Table of Contents}\smallskip

$\hspace{30mm}$ {\footnotesize \parbox[t]{115mm}{
1.{ $\;$} {\bf Introduction}\newline
2.{ $\;$} {\bf Holonomic modules over quantum tori}\newline
3.{ $\;$} {\bf $\hh[W]$-modules}\newline
4.{ $\;$} {\bf Morita equivalence}\newline
5.{ $\;$} {\bf Representations and $G$-bundles on elliptic curves}\newline
6.{ $\;$} {\bf From quantum tori to the Cherednik algebra}\newline
7.{ $\;$} {\bf Operator realization of the Cherednik algebra}\newline
8.{ $\;$} {\bf The Spherical subalgebra}\newline
9.{ $\;$} {\bf Appendix:}\newline
{ $\hphantom{x}\enspace\;\,\,$} {\bf  Springer correspondence 
for disconnected groups}
}}

\vskip 10mm

\section{Introduction.}

\ab We introduce a  non-commutative deformation  of the
algebra  of regular  functions on  a torus.   This  deformation $\hh$,
called \textit{quantum torus algebra},  depends on a complex parameter
$\eps \in \C^*$. 
 We further  introduce a  certain category  $\oo(\hh,\A)$ of
representations of  $\hh$ which are  locally-finite with respect  to a commutative
subalgebra $\A\subset \hh$  whose `size' is one-half of  that of $\hh$
(our  definition  is  modeled   on  the  definition  of  the  category
$\mathcal{O}$ of Bernstein-Gelfand-Gelfand).
  We  classify all simple objects  of $\oo(\hh,\A)$ and
show that any object of $\oo(\hh,\A)$ has finite length.

\ab  In \S3 we  consider quantum  tori arising  from a  pair of
lattices coming from a finite reduced root system. Let $W$ be the Weyl
group of  this root  system. We classify  all simple modules  over the
twisted  group   ring  $\hh[W]$  which  belong   to  $\oo(\hh,\A)$  as
$\hh$-modules.
 In \S4 we  show that  the twisted  group ring  $\hh[W]$ is
Morita equivalent to $\hh^W$, the ring of $W$-invariants.

\ab In  \S5 we establish  a bijection between the  set of simple
modules  over  the  algebra  $\hh[W]$ associated with  a  semisimple
simply-connected group $G$,  and the set of pairs  $(P, \alpha)$, where
$P$  is  a  semistable  principal  $G$-bundle on  the  elliptic  curve
$\mathcal{E}   =   \C^*\!/\eps^{\Z}$,   and   $\alpha$   is   a   certain
`admissible representation' (cf.  Definition 5.4) of the finite
group $Aut(P)/
(Aut\ P)^{\circ}$.

\ab In \S6    we construct a family of algebras $\hdd_{v}$
  depending on a parameter $v$, such that $\hdd_v \simeq \hh[W]$, when
  $v=0$,  and $\hdd_v$ is the double-affine Hecke  algebra,  when
$v=1$. An analogue of Deligne-Langlands-Lusztig conjecture
for the double-affine Hecke algebra $\hdd$ is proposed.
In \S7 we
give an explicit realization of the double-affine
Hecke  algebra  $\hdd$ as  a  subalgebra  of  the twisted  group  ring
$\hh_{_{{\sf   frac}}}[W]$,  for  a  certain   enlargement  $\hh_{_{{\sf
      frac}}}$ of $\hh$. We finally introduce 
  an important \textit{spherical} subalgebra  in $\hdd_v$,
that specializes, at $v=0$, to the subalgebra in $\hh[W]$,
formed by $W$-invariants in the quantum torus.

\vskip 10mm

\section{Holonomic modules over quantum tori.}

Choose to  a finite  rank abelian  group $\V$, referred  to as  a {\it
  lattice},   a   positive   integer   $n$  and   a   skew   symmetric
$\frac{1}{n}\Z$-valued  bilinear  form  $\omega  : \V  \times  \V  \to
\frac{1}{n}\Z$ (where  $\frac{1}{n} \Z$ is  the group of  all rational
numbers of the  form $\frac{a}{n}$, $a \in \Z$).   Associated to these
data is the Heisenberg central extension
\[ 0\to \frac{1}{n}\Z\to {\tilde \V} \to\V\to 0.\]
Here ${\tilde \V}=\V\oplus\frac{1}{n}\Z$ as a set, and the group law on
${\tilde \V}$ is given by 
\[ (v_1, z_1)\circ (v_2, z_2) =(v_1+v_2, z_1+z_2+\omega(v_1,v_2))
\quad, \quad v_i\in \V,z_i\in\frac{1}{n}\Z.\]
Let $\C{\tilde \V}$ denote the group algebra of ${\tilde \V}$ 
formed by all $\C$-linear combinations $\sum_{g\in {\tilde \V}}\, c_g [g]$.
Given a complex number $\eps\in\C^*$ together with a choice of its
$n$-th root $\eps^{\frac{1}{n}}$, we define a {\it quantum torus},
$\hh_{\eps}(\V,\omega)$, as the quotient of $\C{\tilde \V}$
modulo the two-sided ideal generated by the (central)
element $[(0,\frac{1}{n})]-\eps^{\frac {1}{n}}\cdot [(0,0)]$. We write
$e^v$ for the image of $[(v,0)]\in\C{\tilde \V}$ in 
$\hh_\eps(\V,\omega)$. The elements $\{e^v, v\in \V\}$ form a
$\C$-basis of $\hh_{\eps}(\V,\omega)$, and we have
\[e^{v_1}\cdot e^{v_2}=\eps^{\omega(v_1,v_2)}\cdot e^{v_1+v_2}\quad,\quad
\forall v_1,v_2\in\V.\]

{\bf Remarks.}
 1) Replacing $\omega$ by $n \omega$ and $\eps^{\frac{1}{n}}$ by $\eps$
we can reduce to the situation when $n=1$. However, for the later
applications to double affine Hecke algebras it is convenient to work
with the fractional powers of $\eps$ assuming that all appropriate
roots are fixed.

\ab 2) Note  that we  \emph{do  not}  assume that  $\omega$  is a  perfect
pairing,  i.e.,  the map: $\,\V \to \frac{1}{n}Hom (\V,  \Z)\,,$
given by: $\,v\mapsto \omega(v,-)\,$ is injective but is not necessarily
surjective; its image may be a sublattice of finite  index.$\quad\square$
\bigskip

{\bf Lemma 2.1.} {\it If the form $\omega$ is non-degenerate, and
$\eps$ is not a root of unity, then the
algebra $\hh_q(\V,\omega)$ is simple.}
\vspace{.3cm}

{\sl Proof.} Suppose $h= \sum_{i=1}^s  c_ie^{v_i}$ is an element 
of a two-sided ideal $J\sset{\cal H}_\eps(\V,\omega)$, where all
$v_i\in\V$ are distinct, and all the $c_i\in\C$ are
nonzero. We claim that $e^{v_i}\in J$ for every $i$,
whence $J=\hh_\eps(\V,\omega)$ 
since the elements $e^{v_i}$ are invertible.

\ab To prove the claim, we use the non-degeneracy of
$\omega$ and the assumption that all the vectors
$v_i$ are distinct to find an element $v\in \V$
such that $\omega(v,v_i)\not=\omega(v,v_j),$ for
any $i\not= j$. 
Hence, since $q$ is not a root of unity,
we conclude 
$$\eps^{k\cdot\omega(v,v_i)} \not= \eps^{k\cdot\omega(v,v_j)}\quad,
\quad\forall k=1,2,\ldots,\enspace\mbox{\small whenever}\enspace
i\not= j.\eqno(2.1.1)$$
Now, for any $k=0,1,\ldots, $ set $u_k:=
 e^{k\cdot v}h e^{-k\cdot v}\in J$. We have
$$ u_k=
 e^{k\cdot v}h e^{-k\cdot v}= 
\sum c_i\cdot e^{k\cdot v}e^{v_i}e^{-k\cdot v} =
\sum_{i=1}^s c_i\cdot \eps^{k\cdot\omega(v,v_i)}\cdot e^{v_i}.
$$
Observe that the determinant of
the matrix $a_{ik}:=\eps^{k\cdot\omega(v,v_i)}$
is  the Wandermonde determinant
$\prod_{i>j}
(\eps^{k\cdot\omega(v,v_i)}-\eps^{k\cdot\omega(v,v_j)})$.
By (2.1.1) this determinant is non-zero, so that
the matrix is invertible. Hence, each of the elements
$e^{v_1},\ldots,e^{v_s}$ can be expressed as a linear combination
of the $u_0,\ldots,u_{s-1} \in J$, and the claim follows. \sq

{\bf Remark.}
 If  $\eps^m=1$, then the elements $(1-e^{mv}), v\in \V$ are in the 
center of $\hh_\eps(\V,\omega)$, hence any such element 
generates a non-trivial two-sided ideal.

\vspace{.4cm}

\ab  Fix a pair of lattices $\X$, $\Y$ and a non-degenerate pairing 
$\langle \,, \rangle : \X\times\Y\to \frac{1}{n}\Z$.
From now on,  we take  $\V=\X\oplus\Y$, where
   the  form   $\omega$  on
$\X\oplus\Y$ is given by
\[\omega(x\oplus y, x'\oplus y') :=\langle x,y'\rangle -
\langle  x',y\rangle \quad,\quad  x,x' \in  \X,\; y,y'\in  \Y,\] 
 Let $\hh=\hh_\eps(\X\oplus\Y,\omega)$
denote  the corresponding  algebra. The  elements $\{e^x,  x \in\X\}$,
resp. $\{e^y, y\in \Y\}$,  span the commutative subalgebra $\C\X \sset
\hh$, resp.,  $\C\Y \sset  \hh$, and there  is a natural  vector space
(but not {\it algebra}) isomorphism $\,\hh\simeq \C\X\otimes_{_\C}\C\Y
\,.$ The algebra structure is determined by the commutation relations
$$e^y e^x= \eps^{<x,y>} e^x e^y\quad,\quad\forall x\in\X,\,
y\in\Y.\eqno(2.2)$$

\ab We introduce the complex torus
$\;\dis T:= \Hom(\X,\C^*)\;$  so that
$\;\dis \X\simeq$
$  \Hom_{\mbox{{\tiny  alg  group}}}(T,  \C^*)\,.\,$
  Any  element
$x\in\X$ may be viewed as  a $\C^*$-valued regular function $ t\mapsto
x(t)$ on  $T$. For  $y  \in \Y$, the element $n\cdot 
y\in n\cdot\Y$
gives   a well-defined   element   
$\phi_{ny}   \in    \Hom_{\mbox{{\tiny   alg
      group}}}(\C^*, T) = Hom(\X, \Z)$.  We let $\eps^{y} \in T$  be
$\phi_{ny}  (\eps^{\frac{1}{n}})$.  The  assignment  $y\mapsto \eps^y$
identifies  the  lattice  $\Y$  with  a  finitely  generated  discrete
subgroup $\eps^{\Y} \subset T$.

\vspace{.3cm}

\ab Let $A$ be a commutative $\C$-algebra and $\alpha: A\to \C$
an algebra homomorphism, referred to as a {\it weight}.
For an $A$-module $M$, let $M(\alpha):=$
$\{m\in M\,|\,
a m=\alpha(a)\cdot m\,,\,\forall a\in A\}$ denote the corresponding
weight subspace.
\vspace{.5cm}

{\bf Definition.}  Given a $\C$-algebra  $H$ with a  {\it commutative}
subalgebra   $A\sset   H$,  define\hfill\linebreak   $\bullet\enspace$
\parbox[t]{120mm}{$\oo(H,A)$ to be  the category of finitely generated
  $H$-modules $M$ such that the $H$-action on $M$ restricted to $A$ is
  {\it locally finite},  that is for any $m\in  M$ we have $\dim_{_\C}
  A\cdot      m<\infty.\hfill$}      \linebreak      $\bullet\enspace$
\parbox[t]{120mm}{$\os(H,A)$ to be  the full subcategory of $\oo(H,A)$
  consisting of $A$-{\it diagonalizable} $H$-modules, i.e. $H$-modules
  $M$ of the form}
$$M=\bigoplus_{\alpha\in\,\mbox{\tiny Weights of}\,A}\,
M(\alpha)\quad\mbox{\small and}\quad \dim_\C M(\alpha)<\infty\,,\,\forall
\alpha.$$

\bigskip

\noindent
Note that if $A=\C$ then $\oo(H,A)=\os(H,A)$ is just the
category of finitely generated $H$-modules.
\vspace{.3cm}

\ab 
In this section we 
will be concerned with the special case $H=\hh$
$A=\A:=\C\X\sset\hh\,,$
(we also fix $\eps\in\C^*$, not a root of unity).
  Observe that  any object  $M\in  \oo(\hh,\A)$ is
generated  by a  finite dimensional  $\A$-stable subspace.  It follows
that $M$ is finitely generated over the subalgebra $\C\Y\sset\hh$, due
to the vector space factorization $\hh=\C\Y\cdot\C\X$. Since $\C\Y$ is
a Noetherian  algebra, we deduce that, any  $\hh$-submodule $N\sset M$
is  finitely  generated over  $\hh$,  whence $N\in\oo(\hh,\A)$.  Thus,
$\oo(\hh,\A)$ is  an {\it abelian}  category. 
Note the canonical algebra isomorphism
$\,\C\X\,\simeq \,\C[T]\,,$ where $\C[T]$ stands for
the algebra of regular polynomial functions on $T$.
Thus, the
set of weights of the  algebra $\A=\C\X$ is canonically identified with  $T$.

\ab For $\la \in T$, define an $\hh$-module $\RR_{\la}$ as a 
$\C$-vector space with basis $\{v_{\mu}\,,\,
\mu \in \la\cdot \eps^\Y\sset T\}$ and with $\hh$-action given by
$$
e^y(v_{\mu})=v_{\mu\cdot \eps^y}\quad,\quad e^x(v_{\mu})= 
x(\mu)\cdot v_{\mu}.\eqno(2.3)$$
The module $\RR_{\la}$ has the following interpretation. 
Write $I_\mu$ for the maximal ideal in $\C[T]$ corresponding to a
point $\mu\in T$, and let $\C_\mu:=\C[T]/I_\mu$ be the
sky-scraper sheaf at $\mu$. Let $\C[\la\cdot \eps^\Y]:=
\bigoplus_{\mu \in \la\cdot \eps^\Y}\, \C_\mu$
be the (not finitely generated) $\C[T]$-
module formed by all  $\C$-valued, finitely supported
functions on the set $\la\cdot \eps^\Y$. Define an $\hh$-action on
$\C[\la\cdot \eps^\Y]$ by the formulas
$$
e^x (f)  : t \mapsto x(t)\cdot f(t) \quad,  \quad e^y(f): t \mapsto
f(\eps^y\cdot t).\eqno(2.4)  $$
Thus, $x\in  \X$ and $y\in\Y$  act via
multiplication  by   the  function  $x(t)$  and   shift  by  $\eps^y$,
respectively.  It is straightforward  to verify that sending $v_\mu\in
\RR_{\la}\,,\,  \mu  \in  \la\cdot  \eps^\Y $  to  the  characteristic
function of the one-point  set $\{\mu\}$ establishes an isomorphism of
$\hh$-modules  $\,\RR_{\la}\iso  \C[\la\cdot \eps^\Y]\,$  intertwining
the actions (2.3) and (2.4), respectively.

\ab Clearly,  $\RR_\la\in \os(\hh,\A)$.  Moreover, it is  obvious from
the     isomorphism    $\RR_\la\simeq\C[\la\cdot     \eps^\Y]$    that
$\RR_\la\simeq  \RR_\mu$  if  $\mu\in  \la\cdot  \eps^\Y$.  Thus,  the
modules $\RR_\la$ are effectively  parametrized (up to isomorphism) by
the points of the  variety:
$\,\dis\La:=T/\eps^\Y \,.$
When $|\eps| \neq  1$, $\La$ is an abelian  variety.  Observe that the
modules corresponding  to two different points of  $\La$ have disjoint
weights, hence are non-isomorphic.  \vspace{.3cm}

{\bf Proposition 2.5.} (i) {\it $\RR_\la$ is a simple $\hh$-module,
for any $\la\in\La$. Moreover, the set  $\{\RR_\la , \la\in\La\}$
is a complete collection of (the isomorphism classes of)
simple objects of the category $\oo(\hh,\A)$.}

\ab (ii) {\it Any object of the category $\os(\hh,\A)$ is isomorphic
to a finite direct\par
\ab  sum $\bigoplus_{\la\in\La} \,\RR_\la$, in
particular, the category $\os(\hh,\A)$ is semisimple.}

\ab (iii) {\it Any object of the category $\oo(\hh,\A)$
has finite length.}
\vspace{.3cm}

{\sl Proof.} 
Let $M\in \oo(\hh,\A)$. An easy straightforward calculation shows
that, for any non-zero element $m\in M(\la)$, the $\hh$-submodule
in $M$ generated by $m$ is isomorphic to $\RR_{\la}$. This, combined
with the observation preceding the proposition proves part (i).

\ab Since $M$  is finitely generated,  one can find finitely  many weights
$\la_1, \ldots,$
$ \la_s\in T$ such that all weights of $M$ are contained
in $(\la_1\cdot  \eps^\Y) \cup \ldots \cup  (\la_s\cdot \eps^\Y)$ and,
moreover, $\la_i  \neq \la_j\mbox{ mod }\eps^\Y$ whenever  $i\not= j$. 
It follows, since all weights of $M$ are in $(\la_1\cdot \eps^\Y) \cup
\ldots  \cup  (\la_s\cdot \eps^\Y)$,  that  $M$  is  generated by  the
subspace   $\bigoplus_{i=1}^s   M(\la_i)$.   Furthermore,   the   same
calculation as in  the first part implies that  the $\hh$-submodule in
$M$  generated by  this subspace  is isomorphic  to $\bigoplus_{i=1}^s
\RR_{\la_i} \otimes M(\la_i)$. This proves part (ii).

\ab To  prove (iii),  suppose $M\in\oo(\hh,\A)$.  We use
induction on  the minimal dimension $d$ of  an $\A$-invariant subspace
$V  \subset M$ which  generates $M$  over $\hh$.  It follows  from the
definitions  that if  $d=1$  then $M  \simeq  \RR_{\lambda}$ for  some
$\lambda$. If $d>1$,  choose a non-zero 
vector $v \in  V$ of some $\A$-weight
$\lambda$ and note that such  a choice induces a non-zero homomorphism
of  $\hh$-modules  $\RR_{\lambda}  \to  M$. Since  $\RR_{\lambda}$  is
simple,  this  homomorphism  is  necessarily injective.  The  quotient
 $M/\RR_{\lambda}$ is generated  by an  $\A$-invariant subspace
$V/\langle v\rangle$, hence we can apply  the assumption of induction to
this 
$\hh$-module,
and (iii) follows.~\sq \vspace{.5cm}

\pagebreak[3]
\section{$\hh [W]$-modules.}

Let $\Delta \subset \frak{h}$ be a finite reduced root system. Let $W$
be the Weyl group of  $\Delta$ and let $\X \subset \frak{h}^\vee$, $\Y
\subset  \frak{h}$ be a  pair of  $W$-invariant lattices associated
with
$\Delta$, such as e.g.,  the (co)root and weight lattices.  The group $W$
acts  naturally  on $\X$  and  on  $\Y$.  The diagonal  $W$-action  on
$\X\oplus\Y$  makes 
$\,\hh= \hh(\X\oplus\Y)\,$  a left  $W$-module with  $W$-action $w  : h
\mapsto  {^w\!h}\;,\,h\in\hh.$  Write  $\hw$  for  the  subalgebra  of
$W$-invariants.   Further, introduce  a {\it  twisted  group algebra},
$\hh[W]$,  as  the complex  vector  space  $\hh\otimes_\C \C[W]$  with
multiplication:
\[(f\otimes w) \cdot (g\otimes y) = (f\cdot {^w\!\!g})\otimes (w\cdot
y)\quad f,g \in \hh,\; w,y \in W\]
We use similar notation $\hh[W']$ for any subgroup $W'\sset W$,
and view $\C\X$, resp. $\C\Y$, as a commutative subalgebra of
$\hh[W']$ via the composition of imbeddings $\C\X\into\hh\into\hh[W']$.

\ab The group $W$ acts naturally on $T$ and on $\La=T/\eps^\Y$.
Given $\la\in T$, consider its image in $\La$,
and let  $W^{\la} \subset W$ denote the isotropy group of the image of
$\la$. The $W^{\la}$-action on $T$ keeps the subset
$\la\cdot \eps^\Y$ stable, hence we may define $W^{\la}$-action on 
$\RR_{\la}$ by the assignment $w: v_{\mu}\mapsto v_{w(\mu)}$.
This way we make the twisted group algebra, $\hh[ W^{\la}]$, 
act on $\RR_{\la}$.
  
\vspace{.5cm}

  {\bf Theorem 3.1} (cf. \cite[2.1]{LS}).
{\it If $M\in \os(\hh,\A)$, then the restriction of
$M$ to  $\hw$-module is
semisimple, i.e., $M\in \os(\hh^W,\A^W)$.
Furthermore, $\ind_{_\hh}^{^{\hh[W]}} M\in \os(\hh[W],\A)$.}

\vspace{.3cm}
  {\sl Proof.} This follows from Proposition 2.5 and the twisted version of
Maschke Theorem, see  [M, Theorems 0.1 and 7.6(iv)].\sq
\vspace{.5cm}

\ab Let ${\cal M}_{\la}(\hh[W^{\la}],  \A)$ be the full subcategory of
$\os(\hh[W^{\la}], \A)$  formed by the  modules $M$ such that  all the
weights of the $\A$-action belong to the coset $\la \cdot \eps^\Y$.

\ab  Let  $M \in  {\cal  M}_{\la}(\hh[W^{\la}],  \A)$.  Note that  the
subgroup  $W^\la$ does not  necessarly map  the weight  space $M(\la)$
into itself:  if $w  \in W^{\la}$ then  by definition of  $W^{\la}$ we
have  $w(\la)  \in  \la\cdot  \eps^\Y$.  Thus,  it  is  possible  that
$w(\la)\not=\la$ so  that , for  $m\in M(\la)$, the element  $w(m)$ is
pushed out of the $M(\la)$.   We define a "corrected" {\it dot-action}
$w:  m\mapsto w\cdot  m$  of the  group  $W^\la$ on  the vector  space
$M(\la)$ as  follows. As we  have seen by  definition, for any  $w \in
W^{\la}$,  there  exists  a  uniquely determined  $y\in\Y$  such  that
$w(\la)=\la\cdot  \eps^y$. Then, for  $m \in  M(\la)$, put  $w\cdot m=
e^{-y}w(m)$. Here $w(m)\in  M$ stands for the result  of $w$-action on
$m$, and  we claim that  the element $e^{-y}w(m)$ belongs  to $M(\la)$
(while $w(m)$ is not, in general).

\ab Write  $\oo(W^{\la})$ for the  category of finite  dimensional $\C
W^{\la}$-modules.  With the  dot-action of $W^{\la}$ introduced above,
we  may  now  define  a  functor  (cf.  \cite[2.2]{LS})  $\Phi:  {\cal
  M}_{\la}(\hh[W^{\la}], \A) \sto \oo(W^{\la})\,$ by the assignment $M
\mapsto M(\la)$.   On the  other hand, given  a representation  $N$ of
$W^{\la}$  one has an  obvious $\hh  [ W^{\la}]$-action  on $\RR_{\la}
\otimes_{_\C} N$  and this gives a functor  $\Psi: \oo(W^{\la})\, \sto
{\cal M}_{\la}(\hh[W^{\la}], \A)$.

\vspace{.5cm}

{\bf Theorem 3.2.} {\it The functors
$\Psi$ and $\Phi$ are mutually inverse equivalences.}
\vspace{.3cm}
  
{\sl Proof.}  One has $  \Phi \Psi(N) \simeq  N$. If $M$ is  in ${\cal
  M}_{\la}$  and $M(\la)=  \Phi(M)$  then by  theorem  3.1, $M  \simeq
(\RR_{\la})^{\oplus m}$ as $\hh$-module and hence $M=\hh\cdot M(\la)$.
Thus, there  is a morphism of $\hh$-modules  $\psi: \Psi(M(\la))=
\RR_{\la}\otimes_{_\C}  M(\la)  \to  M$  given by  $hv_{\la}\otimes  m
\mapsto h(m)$. The map $\psi$ is injective since $\RR_{\la}$ is simple
over  $\hh$.  One  can  easily   check  that  $\psi$  is  actually  an
isomorphism of $\hh [ W^{\la}]$-modules.\sq

\ab Since $\hh$ is a subalgebra of $\hh[W^{\la}]$ one may
regard $\hh[W^{\la}]$  as a  {\it right} $\hh$-module.  Let ${\widehat
  {W}}^{\la}$  denote  the  set   of  isomorphism  classes  of  simple
$W^{\la}$-modules.   \vspace{.5cm}

{\bf Proposition 3.3} (cf. \cite[2.4]{LS}).
{\it There is  an $\hh[ W^{\la}]$-module decomposition}
$$
 \hh[W^{\la}]\otimes_\hh \RR_{\la} \cong \bigoplus\nolimits_{\chi\in
   {\widehat W}^\la}\; (\RR_{\la}\otimes_{_\C} \chi)^{\oplus d_{\chi}}
\quad,\qquad d_\chi:=\dim\chi
$$
{\it Furthermore, the 
$\hh[ W^{\la}]$-modules $\{\RR_{\la}\otimes_{_\C} \chi\,,\,
\chi \in {\widehat W}^\la\}$ are simple and
  pairwise non-isomorphic.}
\vspace{.3cm}
 
{\sl Proof.} $\Phi(\hh[W^{\la}]\otimes_\hh \RR_{\la})$ is the left regular
 representation of $W^{\la}$.\sq
\vspace{.3cm}

\ab For any $\chi \in {\widehat W}^\la$, put
$V_\chi:=\Psi(\chi)=\RR_{\la}\otimes_{_\C}\chi\;\in\;
{\cal M}_{\la}(\hh[W^{\la}], \A)\,.$ Set

$$
Z_{\chi}\;\,:=\;\, \ind^{^{\hh[W]}}_{_{\hh[W^\la]}} V_\chi\;\,=\;\,
\hh[W]\bigotimes\nolimits_{\hh[W^{\la}]}\;V_\chi\;\,\in \os(\hh[W^{\la}],
\A).
$$
\vspace{.5cm}

{\bf Theorem 3.4} (cf. \cite[2.5]{LS}).
{\it  There is  an $\hh[W]$-module isomorphism}
$$
\hh[W]\otimes_\hh  \RR_{\la}   \cong  \bigoplus\nolimits_{\chi\in
  {\widehat W}^\la}\; {Z_{\chi}}^{\oplus d_{\chi}}
$$
{\it Furthermore,  $Z_{\chi}$ are simple
 pairwise non-isomorphic $\hh[W]$-modules.}
\vspace{.3cm}
 
{\sl Proof.} We have an obvious isomorphism:
\[\hh[W]\otimes_\hh
\RR_{\la}\; \cong\;  \hh[W]\otimes_{\hh[W^{\la}]} \hh[W^{\la}] \otimes_\hh
\RR_{\la}.\] The decomposition of the Theorem now follows from 
Proposition 3.3.
To  prove that  $Z_{\chi}$ are  simple  $\hh[W]$-modules
we write an  $\hh[W]$-module
direct sum decomposition:
\[Z_{\chi}
 \cong \bigoplus\nolimits_{j=1}^s\; w_jV_\chi\quad\mbox{\small and}\quad w_jV_\chi \cong
 (\RR_{h_j(\la)})^{\oplus d_{\chi}},\]
 where $w_1=e, \ldots, w_s,$ are
representatives in $W$ of the right  cosets 
 $W/W^{\la}$.  Any simple $\hh$-submodule of
 $Z_{\chi}$ is contained in some $w_j V_{\chi}$.
 
 \ab   By  Theorem  3.1,   the  $\hh[W]$-module   $  \hh[W]\otimes_\hh
 \RR_{\la}$  is  semisimple.  Therefore,  $Z_{\chi}$, being  a  direct
 summand  of a  semisimple module,  is a  semisimple  $\hh[W]$-module. 
 Hence  $Z_{\chi}$ contains  a simple  submodule $M$  with  a non-zero
 projection from $M$ to $w_jV_\chi$. Viewing $M$ as an $\hh$-module we
 see that $M= \bigoplus_j\; (M \cap w_jV_\chi)$. Since $V_{\chi}$ is a
 simple  $\hh[W^{\la}]$-module,  we  have  $V_{\chi}  \subset  M$  and
 therefore $\bigoplus_j w_jV_\chi \subset M$. Hence $Z_{\chi}=M$.

\ab  Finally, any isomorphism $\theta: Z_{\chi} \to Z_{\psi}$ for some
  $\chi \neq \psi$ maps $V_{\chi}$ to $V_{\psi}$ (just view it as a
  morphism of $\hh$-modules). This would contradict Proposition 3.3.\sq

\vspace{.5cm}

{\bf Proposition  3.5.} {\it Any simple $\hh[W]$-module  $M$ such that
  $\C\X^W$-action  on $M$  is   locally   finite    is   isomorphic    to  
  $Z_{\chi}\,,$ for a certain
$\chi   \in  {\widehat   W}_\la,  \la   \in   \La/W$.}  
\vspace{.3cm}
   
{\sl           Proof.}           We          have           $Z_{\chi}=
\ind_{_{\hh[W^{\la}]}}^{^{\hh[W]}}(V_{\chi})$.   By     Schur    lemma
and Frobenius
reciprocity: $\;\dis\Hom(A, \Res B)=\Hom(\ind A, B)\,,\,$
    it    suffices    to    show    that
 $\,\Res_{_{\hh[W^{\la}]}}^{^{\hh[W]}}(M)$ has  a submodule isomorphic
 to  $V_{\chi}$.  But  the latter  follows from  the proof  of Theorem
 3.2.\sq\medskip

\ab Thus, we have reduced classification of simple $\hh[W]$-modules
to the classification of irreducible representations of the finite group
$W^\la$. The latter group is {\it not} a Weyl group, however.
Therefore its representation theory is not classically known
in geometric terms.
In section 5 we will develop an analogue of "Springer theory"
for $W^\la$, relating
 irreducible representations of
$W^\la$ to semistable $G$-bundles on the elliptic curve $\C^*\!/\eps^\Z$.

\vspace{.5cm}

{\bf Remark 3.6.} Note that one has the following
 alternative definition of $Z_{\chi}$:
 \[Z_{\chi}:= \ind^{^{\hh[W]}}_{_{\A[W^\la]}} (\la\otimes\chi)=
\hh[W]\bigotimes\nolimits_{\A[W^{\la}]}\;(\la\otimes\chi) \,,\]
where $\la$ denotes the one-dimensional $\A[W^{\la}]$-module,
in which the group $W^{\la} \subset \A[W^{\la}]$ acts via the dot-action.

\section{Morita equivalence.}

\ab The  algebra $\hh$  may be  viewed either  as an  $(\hh[W],  \hh^W)$ -
bimodule, $\hh^l$, or as an $(\hh^W, \hh[W])$-bimodule, $\hh^r$.

\vspace{.5cm}

  {\bf Proposition 4.1} (cf. \cite[3.1]{LS}). (i) 
{\it $\hh[W]$ and $\hh^W$ are simple rings. These rings are Morita
equivalent via the following functors}:
\[
{\bold F}: \hh[W]\mod  \sto   \hh^W\mod \quad,\quad M 
\mapsto \hh^r\otimes_{\hh[W]}M\]
\[{\bold I}:\, \hh^W\mod  \sto\hh[W]\mod \quad,\quad 
N \mapsto
\hh^l\otimes_{\hh^W}N \]

(ii) {\it There are functorial isomorphisms}:
$\;{\bold F}(M) \cong \Hom_{\hh[W]}(\hh^l, M) \cong M^W$.
\vspace{.3cm}
  
{\sl Proof.} (i) See \cite[Theorems 2.3 and 2.5(a)]{M}.
(ii) Exercise.\sq
\vspace{.3cm}

\ab Similar results hold for $\hh^{W^{\la}}$- and
$\hh[W^{\la}]$-modules, respectively.
 We write  ${\bold F}_{\la}$ and ${\bold I}_{\la}$ for
the corresponding functors.

\ab Since  $\hh^W$ commutes with $W^{\la}$, we  may regard $\RR_{\la}$
as     a      left     $\hh^W     \times      W^{\la}$-module.     Let
$L_\chi=\Hom_{W^\la}(\chi^*,  \RR_{\la})$   be  the  $\chi^*$-isotypic
component  of the  $\hh[W^{\la}]$-module $\RR_{\la}$.  Notice  that by
Proposition 4.1(ii) we have
\[L_\chi= (\RR_{\la} \otimes_{_\C} \chi)^{W^\la}=
{\bold  F}_{\la}(V_{\chi})=  \hh\otimes_{\hh[W^{\la}]}V_{\chi}= {\bold
  F}(Z_{\chi}).\] Since $\RR_{\la} \cong \bigoplus_{\chi \in {\widehat
    W}_\la}     L_\chi\otimes      V^*_{\chi}$     as     $\hh^W\times
W^{\la}$-modules, we deduce an $\hh[W]$-module decomposition:
$$
  \RR_{\la} \cong \hh\otimes_{\hh[W]}\hh[W]\otimes_\hh \RR_{\la} \cong
 \bigoplus\nolimits_{\chi}\;(\hh\otimes_{\hh[W]}Z_{\chi})^{\oplus d_{\chi}} =
 \bigoplus\nolimits_{\chi}\;L_\chi^{\;\oplus d_{\chi}}.\eqno(4.2)
$$

\vspace{.2cm}

 {\bf Theorem 4.3}  (cf. \cite[3.4]{LS}).
(i) {\it The $\hh^W$-modules $\{L_\chi\,,\,\chi \in
    {\widehat W}^\la\}$ are simple and pairwise
 non-isomorphic.}

\ab (ii) {\it Every simple object of $\oo(\hw,\A^W)$
    is  isomorphic  to  $L_\chi$,  for some  $\chi  \in  {\widehat
    W}^\la$.}  \vspace{.2cm}
  
{\sl Proof.}
(i) Follows from Theorem 3.4 and Morita equivalence.
(ii) Follows from Proposition 3.5 and Morita equvalence.\sq
\vspace{.5cm}

  {\bf Proposition 4.4} (cf. \cite[3.6]{LS}).
{\it If $\RR_{\la}$ and $\RR_{\mu}$ have a simple
 $\hh^W$-submodule in common then $\mu \in W \cdot {\la}$, in which case
$\RR_{\la} \cong 
\RR_{\mu}$.}
 \vspace{.3cm}
  
{\sl Proof.} By Morita equivalence and the identity $\RR_{\la} \cong 
\hh\otimes_{\hh[W]}\hh[W]\otimes_\hh \RR_{\la}$ it's enough to consider the 
$\hh[W]$-modules $\hh[W]\otimes_\hh \RR_{\la}$ and $\hh[W]\otimes_\hh \RR_{\mu}$.
Now consider these modules as $\hh$-modules and apply the 
decomposition $Z_{\chi} \cong \bigoplus_{j=1}^s w_jV_\chi$ from the proof 
of Theorem 3.4. \sq

\section{Representations  and $G$-bundles on elliptic curves}

\ab In   this   section  we   fix $G$, 
a   connected   and
simply-connected complex  semisimple group.  We write $\TT$ for the
{\it abstract} Cartan subgroup of $G$, that is: $\TT := B/[B,B]$,
for an arbitrary Borel subgroup $B\subset G$, see [CG, ch.3].
Let $\W$ denote  the
{\it abstract} Weyl group, the  group acting on $\TT$ and
generated by the given set of simple reflections.
We also fix $\eps\in\C^*$ such that $|\eps| <  1$,
and set $\E=\C^*\!/\eps^\Z$. 

\ab  For any complex reductive group $H$ we let
 $\MM(\E,H)$ denote the moduli space 
of topologically
trivial semistable  $H$-bundles on $\E$.\medskip

{\bf Definition 5.1.}$\;$
A $G$-bundle $P\in \MM(\E,G)$ is called
 `semisimple' if any of the following 3 equivalent conditions hold:

\ab  \vi The structure group of $P$  can be
reduced from $G$ to a maximal torus $T\subset G$;

\ab  \vii The automorphism group $Aut P$ is reductive;

\ab  \viii The {\it substack}
corresponding to the isomorphism class of $P$ is closed
in the {\it  stack} of all $G$-bunles on $\E$.\bigskip

\ab We write $\MM(\E,G)^{ss}$ for the subspace in $\MM(\E,G)$
formed by semisimple $G$-bundles. To each $G$-bundle $P\in \MM(\E,G)$
one can assign its {\it semisimplification}, $P^s\in \MM(\E,G)^{ss}$.
By definition, $P^s$ corresponds to the unique 
closed isomorphism class in the stack  of  $G$-bunles on $\E$
which is contained in the closure of the  isomorphism class of $P$.
This gives the  semisimplification morphism
$ss: \MM(\E,G)\to \MM(\E,G)^{ss}$.
It is known further that there are natural isomorphisms
of algebraic varieties: 
$$\MM^\circ(\E,\TT) \simeq X_*(\TT)\otimes_\Z \E\quad\mbox{\textsl {and}}\quad
\MM(\E,G)^{ss}\simeq \bigl(X_*(\TT)\otimes_\Z \E\bigr)/W\,,\eqno(5.2)$$
where $\MM^\circ(\E,\TT)$ stands for the connected component of the
trivial
representation in $\MM(\E,\TT)$. Moreover, the connected components of 
$\MM(\E,\TT)$ are labelled by the lattice $X_*(\TT)$, and are all
isomorphic to each other.

\ab
By a $B$-structure on a $G$-bundle $P$ we mean a reduction
of its structure group  from $G$ to a Borel subgroup of $G$.
Let $\bb(\E,G)$ denote the moduli space of pairs:
$\{G\mbox{\it-bundle }\, P\in  \MM(\E,G)\,,\,
B\mbox{\it-structure on }\, P\}.\,$ Forgetting the  $B$-structure
gives a canonical morphism $\pi: \bb(\E,G)\too \MM(\E,G)$.
On the other hand, given a  $B$-structure on  $P$ one gets
 a $B$-bundle $P_B$, and push-out via the homomorphism:
$B\onto B/[B,B]=\TT$ gives a $\TT$-bundle on $\E$.
Thus, there is a well-defined morphism of algebraic varieties
$\nu: \bb(\E,G)\too \MM(\E,\TT)$.
Further, set: $\widetilde{G} =
\,\{(x, B)\;\;|\;\; B\;\;\mbox{\it is Borel subalgebra in }\,G\;,\; x\in
B\},\,$ and let $\pi: \widetilde{G}\to G$ be the first projection.

\ab We have the following two commutative diagrams, where the one
on the left is the Grothendieck-Springer "universal resolution"
diagram, cf. e.g. [CG, ch.3], and  the one
on the right is its `analogue' for bundles on the elliptic curve
$\E$:
$$\scriptsize{
\diagram
&\widetilde{G}\dlto_{\pi}\drto^{\nu}&&
  &\bb(\E,G)\dlto_{\pi}\drto^{\nu}&\\
G\drto&&\!\TT\dlto&\!\!
\MM(\E,G)\drto^{ss}&&\!\!\!\!\MM(\E,\TT)\;\;\hphantom{x}\dlto\\
&{\mathtt{Spec}(\C[G]^G)}\simeq\TT/\W&&
& \MM(\E,G)^{ss}\simeq\MM(\E,\TT)/\W&
\enddiagram}
\eqno(5.3)
$$

\ab  Observe that, for any $P \in \MM(\E,G)$ the group $Aut \ P$ acts
naturally
on the set $\bb(\E,G)_{_P} := \pi^{-1}(P)$ of all $B$-structures on $P$.
This induces an action of $Aut \ P/Aut^\circ P$, the (finite) group of connected
components, on the complex top homology
group: $H_{_{top}}(\bb(\E,G)_{_P}, \C)$.
\medskip

\textbf{Definition  5.4} $\,\,$  An  irreducible  representation   of   the group 
$Aut \ P/Aut^\circ P$  is called  `admissible' if  it
occurs in $H_{_{top}}(\bb(\E,G)_{_P}, \C)$ with non-zero multiplicity.
\medskip

\ab One of the main results of this paper is the following

\bigskip
\noindent
\textbf{Theorem 5.5} {\it There exists a bijection between the set of
(isomorphism classes of)
simple objects of $\oo(\hh[W], \A)$ and the set of (isomorphism classes
of) pairs $(P, \alpha)$, where $P\in \MM(\E,G)$,
 and $\alpha$ is an {\rm admissible} representation 
of the  group $Aut \ P/(Aut \ P)^{\circ}$.}

\ab The rest of this section is devoted to the proof of the Theorem.
As a first approximation, recall Proposition 2.5, saying that
simple objects of the category $\oo(\hh,\A)$ are in one-to-one
correspondence with the points of the abelian variety 
$\La=T/\eps^\Z$ which is, by (5.2), nothing but $\MM^\circ(\E,\TT)$.
In the same spirit, it turns out that replacing
algebra $\hh$ by $\hh[W]$  leads to the replacement
of  $\MM(\E,\TT)$ to  $\MM(\E,G)$, as a parameter space for 
simple modules. Specifically, the transition from Proposition 3.5
 to $G$-bundles will be carried out in two steps. In the first step,
we reinterpret the data involved in Proposition 3.5 in terms of loop groups,
and in the second step we pass from loop groups to
$G$-bundles.

\ab
We need some notation regarding formal loop groups.
Let  $\C((z))$, $\C[[z]]$,  $\C[z]$  be the  field  of formal  Laurent
series, the ring of formal  Taylor series and the ring of polynomials,
respectively.    Let    $G((z))$   be   the  group   of   all
$\C((z))$-rational points of $G$,  and similarly for $G[[z]]$, $G[z]$. 
We consider $\eps$-conjugacy classes in $G((z))$, i.e.  
$G((z))$-orbits on itself under $\eps$-conjugation:
$\;g(z) : \,h(z)\mapsto g(\eps z)  h(z) g(z)^{-1}\,.$
A   $\eps$-conjugacy class,
is said to be
\emph{integral} if it contains at  least  one  element  in $G[[z]]$.

\ab Fix
a  Borel subgroup $B=T\cdot U \subset  G$, where $T$ is a maximal torus
of $G$ and
$U$ is the unipotent radical of $B$. By [BG, Lemma  2.2] we have:

\medskip
{\bf Jordan $\eps$-normal  form for  ${\bold G[[z]]}$.}$\;$ {\it
Any   element  $h  \in   G[[z]]$  is
$\eps$-conjugate to a product $s \cdot b(z)$, where
$s \in T$ is a constant loop, and $b \in U[z]$ are such that:

 \ab ${\bold{(J1)}}\enspace$ $b(\eps z) \cdot s  = s \cdot b(z)$,

 \ab ${\bold{(J2)}}\enspace$ $\ad{s} (v) = \eps^m v\,,$
for some $v \in \Lie G\,,\,m > 0\enspace\;\Longrightarrow\enspace\;
v \in \Lie{U}$.}
 \bigskip

\ab For any group $M$, we  write $M^\circ$ for the identity connected 
component of $M$, and $Z_M(x)$ for the  centralizer
of an element $x$ in $M$.
Given $h \in G((z))$ we write $G_{\eps, h}$ for the
$\eps$-centralizer of $h(z)$ in $G((z))$:
\[
G_{\eps, h}\, := \, \{g(z) \in G((z)) \quad | \quad g(\eps z)
h(z) g(z)^{-1} = h(z)\}
\]
Let $W_G=N_G(T)/T$ be the Weyl group of $(G,T)$. Given
$s\in T$, write $\la(s)$ for its image in $\La=T/\eps^\Z$,
and let $W^{\la(s)}$ denote the isotropy group of the
point $\la(s)\in T/\eps^\Z$ under the natural $W$-action.
\medskip

{\bf Theorem 5.6.} {\it Let  $h=s\cdot b\in G[[z]]$ be
 written in its $\eps$-normal  form. Then we have:
$\;G_{\eps, h} = Z_{_{G_{\eps, s}}}(b)\;.$ Furthermore,

\ab \vi $G_{\eps, s}$ is a finite-dimensional reductive group
isomorphic to a (not necessarily connected)
    subgroup $H\subset G$ containing the maximal torus $T$.

\ab \vii There exists a  unipotent element
     $u\in H$, uniquely determined up to conjugacy in $H$,
     such that under the isomorphism in \vi we have:
     $\,G_{\eps, h} = Z_{_{G_{\eps, s}}}(b) \iso Z_H(u)\;.$

\ab \viii The group $W^{\la(s)}$ is isomorphic to
$W_H :=  N_H(T)/T$, the `Weyl group' of the disconnected group
$H$.}
\medskip

The proof of the Theorem will follow from Lemma 5.11 and Proposition 5.13
given later in this section.\medskip

{\bf From loop group to $\bold G$-bundles.}$\;$ In [BG] we have
constructed a  bijection:
$$\MM(\E,G)\enspace\stackrel{\Theta}{\longleftrightarrow}\enspace
\mbox{\it integral}\;\;\eps\mbox{\it-conjugacy classes in }\;
G((z))\;.\eqno(5.7)$$

\ab Let $P=\Theta(h)$ be the $G$-bundle corresponding to
a $\eps$-conjugacy class of $h\in G((z))$,
and $P^s=ss(P)$ its semisimplification. 
Without loss of
generality we may assume that $h$ is written in
its $\eps$-normal  form: $h=s\cdot b$. Using Theorem 5.6
it is easy to verify that under the bijection (5.7) we have:

\ab $\bullet\quad P^s= \Theta(s)\quad\mbox{and}\quad
    Aut \ P^s\, \simeq\, G_{\eps,s}\, \simeq\, H\,\subset G\,.\hfill(5.8.1)$

\ab $\bullet\quad Aut \ P\, \simeq\, G_{\eps,h}\,\simeq\, Z_H(u)\,.$\hfill(5.8.2)\break

Further, recall the variety $\bb(\E,G)_{_P}$ of all $B$-structures
on $P$, see (5.3). Let $\bb(\E,G)_{_P}^\circ$ denote a connected component
of $\bb(\E,G)_{_P}$. Write
 $\bb(H)$ for the Flag variety of the group $H$,
and $\bb(H)_u$ for the Springer fiber over $u$, the $u$-fixed point set
in $\bb(H)$. Then we have:

\ab $\bullet\quad \bb(\E,G)_{_P}^\circ \simeq \bb(H)_u\,.\hfill(5.8.3)$

Furthermore, the natural $Z_{H^\circ}(u)$-action on
$\bb(H)_u$ goes under the isomorphism above and the
imbedding: $\,Z_{H^\circ}(u)
\into Z_H(u)= Aut \ P\,$ to the
natural $Aut\  P$-action on ${\mathcal{B}}(P)$.\medskip

\ab By isomorphism (5.8.3), one identifies
the action of the finite group $Z_u(H^\circ)/Z^\circ_u(H)$ 
on:  $\,H_{_{top}}({\mathcal{B}}_u,\C),$  the  top homology,
with the 
  action of the corresponding subgroup of $Aut\ P/Aut^\circ P$ 
on:  $\,H_{_{top}}(\bb(\E,G)_{_P}^\circ)$.
It follows that an
  irreducible  representation   of  $Aut\ P/Aut^\circ P$ is
admissible in the sense of Definition 5.2 if and only if
the restriction of the corresponding  representation  of
$Z_u(H)/
Z_u^{\circ}(H)$  to the
subgroup $Z_u(H^\circ)/Z^\circ_u(H) \subset Z_u(H)/
Z_u^{\circ}(H)$  is isomorphic to  a direct sum
of irreducible representations which have non-zero multiplicity in the
$Z_u(H^\circ)/Z_u^{\circ}(H)$-module $H({\mathcal{B}}_u)$.

\ab Finally, we observe that
the  isotropy group $W^{\la(s)}$ occurring in part \viii of Theorem 5.6
is exactly the group whose irreducible representations  label
the simple objects of the category $\oo(\hh[W], \A)$, see
Proposition 3.5. Thus, according
to the isomorphism $W^{\la(s)}\simeq W_H$ of Theorem 5.6(iii),
we are interested in a parametrisation of irreducible representations
of the   group $W_H$. Such a parametrisation is provided by
a version of the  Springer
correspondence for disconnected reductive groups,
developed in the
last section (Appendix) of this paper. 
This concludes an outline of the proof of Theorem 5.6.
\bigskip

\ab We now begin a detailed exposition, and
recall  the Bruhat  decomposition  for the group $G[z,  z^{-1}]$.   Let  $G_1[z]
\subset G[z]$ denote the subgroup of loops equal to $e \in G$ at $z=0$
and  denote  by  $\mathcal{U}^+$  the  subgroup $U  \cdot  G_1  [z]$.  
Similarly, $\mathcal{U}^-$ will denote $U^- \cdot G_1[z^{-1}]$ where
$U^- \subset G$ is the unipotent subgroup opposite to $U$ and
$G_1[z^{-1}]$ is the kernel of evaluation map $G[z^{-1}] \to G$ at $z
= \infty$.

\bigskip
\noindent 
\textbf{Proposition 5.9.} \emph{(cf.  [PS, Chapter  8]) Any  element of
  $g(z) \in G[z, z^{-1}]$ admits a \emph{unique} representation of the
  form
$$
g(z) = u_1(z) \cdot \lambda(z) \cdot n_w \cdot  t \cdot u_2(z)
$$
where  $u_1  (z),   u_2(z)  \in  \mathcal{U}^+$,  $\lambda(z)  \in
{\Y} =  Hom_{_{alg}} (\C^*, T)$, $t  \in T$, $w \in  W$ and $u_2
(z)$,  in addition,  satisfies  $[\lambda(z) n_w]  \cdot u_2(z)  \cdot
[\lambda(z) n_w]^{-1} \in \mathcal{U}^-$.} \hfill $\square$

\smallskip
\noindent
\textbf{Corollary 5.10.} \emph{The $\eps$-conjugacy classes that
intersect 
$T\subset G((z))$ are parametrized by ${\mathbf{\Lambda}}/W$.}

\bigskip 
\noindent
\textit{Proof.}
Suppose that $s \in  T$ is $\eps$-conjugate to $s' \in T$
by an element  $g (z) \in G((z))$. Rewriting this  in the form $g(\eps
z) s =  s' g(z)$, then using the  above decomposition and its
uniqueness,  we  obtain  $s'  = w(s)  \cdot  \lambda(\eps)$.  
Conversely, for any $w \in W$ and $\lambda \in \mathbf{Y}$, the element
$s$ is conjugate  to $w(s) \cdot \lambda(\eps)$ by  the element $g(z) =
\lambda(z) \cdot n_w$. \hfill $\square$

\medskip
\ab  Uniqueness of the $\eps$-normal Jordan form follows from

\bigskip
\noindent
{\textbf{Lemma 5.11.}} \emph{Suppose that two loops $s \cdot b (z)$ and
  $s'\cdot b'(z)$  satisfy the Jordan form conditions (J1)-(J2),
 and  that $f(\eps z)
  (s \cdot  b (z)) f(z)^{-1}  = s' \cdot  b'(z)$ for
  some $f(z) \in G((z))$.  Then
  $f(\eps z)\cdot  s\cdot  f(z)^{-1}  = s' \quad \textrm{and} \quad f(z)
  b(z) f(z)^{-1} = b'(z).$}

\medskip
\noindent 
\textit{Proof.}  Choose  a faithful representation:  $G \to GL(V)$,  and a
basis in  $V$ such that $U$  maps to upper-triangular  matrices and $T$
maps to diagonal  matrices.  We may assume without loss of generality
that the  loops $s \cdot b(z)$ and
 $s' \cdot b'(z)$  are both maped into  upper-triangular
matrices, $A(z)$ and   $A'(z)$, resp.

\ab First we consider the case when  all diagonal entries of $A(z)$ (resp. 
$A'(z)$) differ only by powers of $\eps$, i.e.  when they are
of the  form $a  \eps^{m_1}, \ldots, a  \eps^{m_k}$, where $k$  is the
dimension of $V$  and $m_1 \geq m_2 \geq \ldots  \geq m_k$,
due to  Jordan form condition (J2).
Further, by the  Jordan form condition (J1)  all entries $A_{ij}$
above  
the diagonal are of
the form $\alpha_{ij}  z^{m_i - m_j}$, $i < j$,  $\alpha_{ij} \in \C$. 
Let also  $a' \eps^{n_1}, \ldots,  a' \eps^{n_k}$ be
the  diagonal entries  of $A'(z)$  and $F(z)$  be  the matrix
correspoding to $f(z)$.

\ab We prove by  descending induction on $i-j$ that $F_{i,j}  = c z^l,$ for
an appropriate  constant $c$  and an integer  $l$, depending  on  $i,j$. Our
proof is  based on the simple  observation that, for  any constant $B$
and any  integer $l$, the equation $x(\eps  z) = \eps^l x(z)  + B z^l$
admits a solution in $\C((z))$ iff $B = 0$, in which case the solution
has to be $x(z) = c z^l$, $c \in \C$.

\ab The largest value  of $i -j$, attained for  $i=k, j=1$, corresponds to
the  lower  left corner  element  $F_{k,1}  (z)$.   From the  equation
$F(\eps z)  A(z) =  A'(z) F(z)$  one has  $F_{k,1} (\eps  z) a
\eps^{m_1} =  F_{k,1} (z) a'  \eps^{n_k}$.  If the  ratio $a/
a'$ is not a power of $\eps$, this equation, as well as other
equations  considered below,  has only  zero solution  (which  gives a
non-invertible matrix $F(z)$).  Hence  we can assume that $a'
=  a  \eps^r$ for  some  integer $r$.   Then  the  above equation  for
$F_{k,1}(z)$ implies that $F_{k,1}(z) =  \phi_{k,1} z^{n_k - m_1 + r}$
for some constant $\phi_{k,1}$.

\ab We  use this  expression  for  $F_{k,1}$ to  write  the equations  for
$F_{k-1, 1}$ and $F_{k, 2}$, then write the equations for $F_{k-2,
  1}$, $F_{k-1, 2}$, $F_{k, 3}$, etc.
In general, by descending induction  on $i-j$ (
ranging from $i-j= k-1$ to $i-j
= -k+1$)  one obtains equations of  the type 
$$F_{i,j} (\eps z)  = \eps^{n_i -
  m_j + r} g_{i,j}  (z) + C z^{n_i - m_j + r}$$
for some constant $C$ depending on $i, j$ and the previously computed
values of $g_{s, t}$.
As before, this  leads to 
$$
C=0 \quad \textrm{and} \quad F_{i,j}(z)  = \phi_{i,j} z^{n_i -
  m_j +  r}, \quad \phi_{ij} \in \C
$$
This equation  implies that $f(\eps z)\cdot  s\cdot  f(z)^{-1}  = s'$ and
$f(z) b(z)  f(z)^{-1} = b'(z)$ is an immediate consequence.

\ab In the general case, by 
 Jordan form condition (J1)   one can choose a basis of $V$ so
that $A(z),  A'(z)$ will have  square blocks as in  the first
part of the proof  (``$\eps$-Jordan blocks'') along the main diagonal,
and zeros everywhere else. We can assume that any two diagonal entries
which differ by a power of $\eps$, belong to the same block.  A direct
computation  shows that,  up to  permutation of  blocks in  $A(z)$ and
$A'(z)$, the conjugating matrix $F(z)$ also has square blocks
along the  main diagonal and zeros  everywhere else. Now  we apply the
above  argument to  each  individual  block to  obtain  the result  in
the general case.  \hfill $\square$

\bigskip
\noindent
\textbf{Corollary 5.12.} 
\emph{\vi The assignment: $s \cdot b(z) \mapsto s$ descends
to a well-defined map}$\;$
\emph{$ \Phi: \;\; \{\textrm{integral } \eps-\textrm{conjugacy
      classes in } G((z)) \} \too {\mathbf{\Lambda}}/W.$}

\medskip
\ab
\emph{\vii  Let $s  \cdot  b(z)$ be  a Jordan $\eps$-normal form,
and $\lambda  \in
  {\mathbf{\Lambda}}/W$       the   image   of   $s    \in   T$   in
  ${\mathbf{\Lambda}}/W$.  Then the set
$\Phi^{-1}(\lambda)$ can be identified
  with those  (ordinary) conjugacy classes  in 
$G_{\eps,  s}$, 
the $\eps$-centralizer of $s$ in $G((z))$,
   which  have  non-trivial intersection  with $U[z]$.}  
\hfill $\square$

\bigskip
{\textbf{Remark.}   We   will  see  below  that  $G_{\eps,   s}$  is  a
finite-dimensional  reductive group  and that  $\Phi^{-1}(\lambda)$ is
nothing but the  set of unipotent conjugacy classes  in this reductive
group.}

\medskip\ab  
Now  we  begin to  study  the  automorphism  group of  the  $G$-bundle
$P^s$ associated  to $s  \in T$.  To  describe $G_{\eps,  s} \simeq 
Aut\  P^s$ first recall  that  by  [BG,  Lemma 2.5],  $G_{\eps,s}$
consists  of polynomial loops,  i.e.  $G_{\eps,s}  \subset G[z,  z^{-1}]$. 
Thus, there is a well-defined  evaluation map $ev_{z=1}: G_{\eps, s}
\to  G$ sending  a polynomial  loop  to its  value at  $z=1$.  Let  $H
\subset G$ be the image of $G_{\eps, s}$. Write $N_H(T)$ for the
normaliser
of $T$ in $H$, and $W_H :=  N_H(T)/T$ for the `Weyl group'
of the (generally {\it diconnected}) group $H$.

\bigskip
\noindent
\textbf{Proposition 5.13.} 
\emph{
\vi The evaluation map $ev_{z=1}: G_{\eps, s} \to G$ is injective;}

\ab
\emph{
\vii The idenity component $H^{\circ}$ of $H$ equals
 the connected reductive subgroup of $G$ corresponding to the root
subsystem $\Delta_{\eps,  s} \subset \Delta$ of  all roots  $\alpha  \in \Delta
\subset Hom(T, \C^*)$ for which $\alpha(s)$
is an integral power of $\eps$;}

\ab
\emph{\viii The  group $W_H$  is isomorphic to
  the subgroup $W^{\lambda} \subset W$ of all $w \in W$ which fix
  $\lambda \in {\mathbf{\Lambda}}=\C^*\!/\eps^\Z$, the image  of $s\in T$.}

\bigskip
\noindent 
\emph{Proof.} 
The equation $g(\eps z) \cdot s \cdot g(z)^{-1} = s$ can be rewritten 
as $g(\eps z) = s \cdot g(z) \cdot s^{-1}$. We decompose $g(z)$ as in
Proposition 5.9 and, using the uniqueness of this decomposition, obtain
$$
u_1 (\eps z) = s u_1(z) s^{-1}, \quad u_2 (\eps z) = s u_2(z) s^{-1},
\quad s = w(s) \cdot \lambda(\eps).
$$
Rewrite  $u_1(z)   \in  \mathcal{U}^+$   as   $\;\displaystyle  exp
\big(\sum\nolimits_{k=0}^\infty g_k  z^k\big),\;$ then $\ad{s} (g_k)  = \eps^k$ due
to the first equation.  In particular, only finitely many of $g_k$ are
non-zero  and by  by Jordan form condition (J2),
  $g_k \in  \Lie(U)$. Hence  $u_1(z) \in
U[z]$ and, since different eigenspaces of $\ad{s}$ on $\Lie(U)$ have zero
intersection,   $u_1(z)$  is   uniquely   determined  by   $u_1(1)  =   exp
(\sum_{k=0}^N g_k)$. The same  argument applies to $u_2(z)$. Moreover,
since  $u_2' =  \big[\lambda(z) n_w  \big] u_2(z)  \big[\lambda(z) n_w
\big]^{-1} \in U^-  \cdot G_1 [z^{-1}]$ and $\ad{s}  u_2'(z) = u_2'(\eps
z)$, we can repeat the argument once more and conclude that $u_2'(1) =
n_w u_2(1) n_w^{-1} \in U^-$.

\ \ \ 
\ab Now we can show that $g(z)$ is determined by $g(1) = u(1) n_w t
u_2(1)$.   In fact,  since  $u_1(1),  u_2(1) \in  U$  and $n_w  u_2(1)
n_w^{-1} \in  U^-$, the usual  Bruhat decomposition for $g(1)  \in G$
implies  that $n_w$,  $u_1(1)$ and  $u_2(1)$ are  uniquely  determined by
$g(1)$, hence  $u_1(z)$ and $u_2(z)$  are uniquely determined by  $g(1)$. 
The element $\lambda(z)$  can be reconstructed from $a$  and $w$ since
$s = w(s) \cdot \lambda(\eps)$ and $\eps$ is not a root of unity. The
proposition follows. \hfill $\square$

\bigskip
\noindent
\textbf{Example.} The following example, showing that the component
group $H/H^{\circ}$ can in fact be nontrivial, was kindly communicated
to us by D. Vogan.

\ab Recall  that for the  root system  of type  $D_4$, the  coroot lattice
$\mathbf{Y}$  can be  identified  with the  subgroup  of the  standard
Eucledian lattice $L_4  = \langle e_1, e_2, e_3,  e_4 \rangle$, $(e_i,
e_j) =  \delta_{ij}$ formed by all  vectors in $L_4$ with  even sum of
coordinates. Then the  set of coroots is identified  with $\pm e_i \pm
e_j$, $i \neq j$, and the  Weyl group $W$ acts by permuting the $e_i$,
and changing the  sign of any even number of the  basis vectors $e_i$. 
The  choice  of  the  simple  coroots $\alpha_1^\vee  =  e_1  -  e_2$,
$\alpha_2^\vee   =  e_2  -   e_3$,  $\alpha_3^\vee   =  e_3   -  e_4$,
$\alpha_4^\vee   =   e_3   +   e_4$   identifies   $T   =   \mathbf{Y}
\otimes_{\mathbb Z} \C^*$ with $(\C^*)^4$. Now consider the element $s
=  (-1, \sqrt{\eps},  -1, -  \sqrt{\eps}) \in  (\C^*)^4 \simeq  T$.  A
straightforward calculation shows that,  in the notations of the above
proposition, $W_H = \{ \pm 1\}$ while $\Delta_{\eps, s}$ is empty.

\bigskip
\noindent
\textbf{End of proof of Theorem 5.5.}\quad
By Proposition 3.5 we have to establish the correspondence between the
set of pairs $(\lambda, \chi)$ where $\lambda \in {\mathbf{\Lambda}}$,
$\chi \in  \widehat{W}^{\lambda}$, and the set of  pairs $(P, \alpha)$
as in the statement of Theorem 5.5. Take any lift $s \in T$ of
the  element   $\lambda  \in  {\mathbf{\Lambda}}$   and  consider  the
$G$-bundle   $P^s$  corresponding   to  $s$,   together   with  its
automorphism  group  $H$.   By  Proposition  5.13  (iii)  and  Springer
Correspondence  (see Appendix),  the representation  $\chi$  defines a
unipotent  orbit  in  $H$-orbit   in  $H^{\circ}$  together  with  the
admissible representation $\alpha$ of  the centralizer of any point
$u$ in this  orbit. The element $u \in  H$ corresponds via Proposition
5.13 (i)  to a certain  loop $b(z)$, such  that $s \cdot b(z)$  is a
$\eps$-normal  form. The   bundle  $P$   corresponds  to   the
$\eps$-conjugacy class of $s \cdot b(z)$. \hfill $\square$

\vspace{.3cm}

\ab It will be convenient for us in the next section to reinterpret
the parameters $(P, \alpha)$ entering
Theorem 5.5 in a different way as follows.
First, giving $P\in \MM(\E,G)$ is equivalent, according to (5.7), to giving
the $\eps$-conjugacy class of an element $h(z)\in G[[z]]$.
Using the Jordan $\eps$-normal form, write
$h(z)=s\cdot b(z)$, where $s\in T$, is a semisimple element in
$G$, the subgroup of constant loops. Furthermore, by Theorem 5.6
we have: $Aut \ P/Aut^{\circ}P =
Z_{G_{\eps,s}}(b)/Z^{\circ}_{G_{\eps,s}}(b)$.

\ab  Let $Q=ss(P)$ be the semisimplification of $P$.
By (5.8.1), this is the $G$-bundle on $\E$ corresponding,
under the bijection (5.7), to the constant
loop $s$. Let $G_{_Q}$ denote the associated vector bundle
on $\E$ corresponding to the principal $G$-bundle
$Q$ and the adjoint representation of the group $G$.
By construction, $b(z)$ is a polynomial loop
with unipotent values that $\eps$-commutes with $s$.
Hence $b(z)$  gives
rise to a unipotent automorphism $\hat{b} \in Aut \ Q$.
 This way one obtains a  bijection:
$$\MM(\E,G)\enspace\longleftrightarrow\enspace
\Big\{
{{{{\mbox{\footnotesize{semisimple 
$G$-bundle $Q\in\MM(\E,G)^{ss} $}}}}}\atop{{\mbox{
\footnotesize{and a unipotent element $u\in Aut \ Q$}}}}}
\Big\}\,.\eqno(5.14)
$$
It is not difficult to show that the set $\bb(\E,G)_{_P}$, see (5.8.3)
 gets identified,
under the bijection above, with the set of $u$-stable
$B$-structures on the
$G$-bundle $ss(P)$.
\bigskip

\ab 
Fix $\eps\in \C^*$, which is not a root of unity.
An element of the group $G((z))$ will be called
{\it $\eps$-semisimple}, resp. {\it $\eps$-unipotent},
 if it is $\eps$-conjugate to a constant
semisimple loop, resp. conjugate (in the ordinary sense)
to  an 
element of $U[z]$. Write $G((z))^{^{\eps{\sf{-ss}}}}$ and
$G((z))^{^{\eps{\sf{-uni}}}}$ for the 
sets of $\eps$-semisimple and $\eps$-unipotent elements, respectively.
Given $h(z)\in G((z))$, recall the notation
$G_{\eps,h}$ for the $\eps$-centralizer of $h$ in $G((z))$,
 and for any $u(z)\in G((z))$, put
\[Z_{\eps,h}(u)=\{g(z)\in   G((z))\;\;\Big|\;\;
g(\eps z)h(z)=h(z)g(z)\;\;\&\;\;g(z)u(z)=u(z)g(z)\}\]
a simultaneous `centralizer' of $h(z)$ and $u(z)$.
If $h$ is $\eps$-semisimple and $u$ $\eps$-commutes with $h$, then
the group $Z_{\eps,h}(u)$  acts on $\bb(G_{\eps,h})_u$, the $u$-fixed
point
set in the Flag variety of the finite-dimensional reductive
group $G_{\eps,h}$, see  Theorem 5.6(i).
This gives a $Z_{\eps,h}(u)/Z^\circ_{\eps,h}(u)$-action 
on 
$H_*(\bb(G_{\eps,h})_u)$, the {\it total} homology.
An irreducible  representation of the
component group $Z_{\eps,h}(u)/Z^\circ_{\eps,h}(u)$
is said to be  {\it admissible} if it occurs in
$H_*(\bb(Z_{\eps,h})_u)$ with non-zero multiplicity.
We let $\widehat{Z_{\eps,h}(u)/Z^\circ_{\eps,h}(u)}$
denote the set 
of admissible
$Z_{\eps,h}(u)/Z^\circ_{\eps,h}(u)$-modules (cf. Definition 5.2
and the paragraph below formula (5.8.3)).

\ab We now consider the following set:
$$
\MMM=\Big\{(s,u,\chi)\enspace\Big|\enspace
{{s\in G((z))^{^{\eps{\sf{-ss}}}}\,,\, 
u\in G((z))^{^{\eps{\sf{-uni}}}}}\atop{
s(z)u(z)s(z)^{-1}= u(\eps z)\,,\,\chi \in 
\widehat{Z_{\eps,s}(u)/Z^\circ_{\eps,s}(u)}}}\Big\}\eqno(5.15)
$$
\ab Thus, we can reformulate Theorem 5.5 as follows\medskip

{\bf Theorem 5.16.} {\it There exists a natural bijection between
the set of isomorphism classes of simple objects 
of $\oo(\hh[W], \A)$ and the set of $\eps$-conjugacy
classes in $\MMM$.}
\bigskip

\section{From quantum tori to the Cherednik algebra}

\ab Let $\Delta\in \X$ be a finite reduced root system, and
$W$ the corresponding finite Weyl group
generated by the simple reflections $s_i\,,\, i=1,\dots, l $.
Let
$\HH$ be the Hecke algebra associated to $W$. Thus, $\HH$ is a
free module  over the Laurent polynomial ring, $\C[\zht,\zht\inv]$,
with the standard $\C[\zht,\zht\inv]$-basis $\{\T_w, w\in W\},$
see [KL]. The base elements
$\T_i\,,\, i=1,\dots, l, $ corresponding to the simple reflections
$s_i\in W$
generate $\HH$, satisfy the braid relations and
the quadratic identity:
$\;\dis
(\T_i-\zht)(\T_i + \zht^{-1})=0\;,\; i=1,\dots, l.\,
$

\ab Write $\X$ and $\Y$ for the root   and co-weight   lattices
of our root system, respectively. 
Set  $T:= \C^*\otimes_{_\Z}\Y$, the corresponding  torus, and
 let  $\C (T)$ denote the algebra  of rational
(meromorphic)
functions
on $T.$ As before, we    write weights $\lambda$  as functions on $T$ using  
the notation $e^\lambda$.

\ab Let $\Raf$ be the affine root system corresponding to
the extended Dynkin diagram of $\Delta$,
and $\dT=T\times\C^*$ the corresponding torus.
We write $\eps$ for the function on $T\times\C^*$ given
by the second projection. Thus, we identify
$\C[\eps,\eps^{-1}]$ with the coordinate ring of the
group $\C^*$,  and will write $q$ for the value of the function
$\eps$ at a particular
point.
Let $\Waff = W\ltimes \Y$ be  the (extended) affine
Weyl group. The affine Hecke algebra $\hd$ associated
to the affine
Weyl group has
$l+1$ generators $\T_i\,,\, i=0,1,\ldots,l$. The algebra $\hd$ 
has a standard faithful representation on $\C [\dT]$ such that
the  operators $\T_i\,,\, i=0,1,\ldots,l,$ are realized by the
Demazure-Lusztig operators:
$$\hT_i = \zht s_i + {{\zht-\zht^{-1}}\over{e^\alpha - 1}}(s_i
- 1)\quad,\quad i=0,1,\ldots,l.\eqno(6.1)$$ 

\ab Recall that the {\it double-affine Hecke algebra}, $\hdd$
introduced by Cherednik, see  [Ch] and also [Ki],
may be defined as the subalgebra of $\C[\zht, \zht^{-1}]$-linear endomorphisms
of $\C(\dT)[\zht, \zht^{-1}]$ generated by multiplication operators 
by the elements of $\C [\dT]$ and by the $l+1$ operatrs (6.1).
This algebra will be referred to as
the `Cherednik algebra', for short.

\ab We now construct a  family of associative algebras $\hddv$,
depending on
a parameter $v\in\C$. Specifically, we let $\hddv$
be the subalgebra of $\C[\zht, \zht^{-1}]$-linear endomorphisms
of $\C(\dT)[\zht, \zht^{-1}]$ generated by multiplication operators 
(by  elements of $\C [\dT]$), and by the following $\,(l+1)\,$ 
$v$-deformed operators:
$$\hT_{i,v} = \zht s_i + {{v(\zht-\zht^{-1})}\over{e^\alpha - 1}}(s_i
- 1)\quad,\quad i=0,1,\ldots,l.\eqno(6.2)$$ 
Then $\;\dis {\hT_{i,v}}^2=v(\zht-\zht^{-1})\hT_{i,v} + v + \zht^2(1-v),\;$
and the $\hT_{i,v}$ satisfy braid relations, by \cite{BE}.
We note that the eigenvalues of $\hT_{i,v}$ are
$\zht$ and $v(\zht-\zht^{-1})-\zht$.
Clearly, for $v=1$  we
have: $\hddv = \hdd$ and, moreover,
$\hddv\simeq\hdd$,  for any $v\neq 0$.
On the other extreme,
for $v=0$ we have:
${\hdd}_{\big|v=0}\simeq \hh[W]$, where $\hh=\hh(\X\oplus\Y)$ is the quantum torus
considered in \S2.
Thus, we can interpret 
the Cherednik algebra as a deformation of the algebra $\hh[W]$.
Note that the usual finite and affine Hecke algebras are
naturally deformed as subalgebras in $\hddv$.\medskip

{\bf Remark.}  In the notation of [GKV], the algebra
$\hddv$ may be characterised as the one associated 
to the vanishing condition
$T_{\alpha,1-v(1-\zht^{-2})}.$$\quad\square$
 \bigskip

\ab We now turn to representation theory of the Cherednik algebra.
We are interested in parametrizing all simple objects of the category
${\mathcal{M}}(\hdd\,,\,\C\X)$.
By a standard argument based on Schur lemma, 
the parameters $\eps$ and $\zht$ specialize to scalars in any simple 
$\hdd$-module from  the category ${\mathcal{M}}(\hdd\,,\,\C\X)$.
Thus, for any $q,t\in \C^*$, we may consider the subcategory
${\mathcal{M}}_{\!_{q,t}}(\hdd\,,\,\C\X)$ of those $\hdd$-modules
on which $\eps$ acts as $q$, and $\zht$ acts as $t$.
Thus, we fix $q\,,\,t$, and  assume from now on that both
$q$ and $t$ are {\sf not roots of unity}.
Notice that a $\epq$-unipotent element
can be rased to the complex power $t$.

\ab We introduce the following deformation of the
set $\MMM$, see (5.15):
$$
\MMM_t=\Big\{(s,u,\chi)\enspace\Big|\enspace
{{s\in G((z))^{^{\epq{\sf{\tiny{-}ss}}}}\,,\, 
u\in G((z))^{^{\epq{\sf{\tiny{-}uni}}}}}\atop{
s(z)u(z)s(z)^{-1}=u(\epq z)^t\,,\,\chi \in 
\widehat{Z_{\epq,s}(u)/Z^\circ_{\epq,s}(u)}}}\Big\}
\eqno(6.3)
$$
\medskip

{\bf Higgs bundle interpretation.}\quad It is instructive to reformulate the data
consisting of a triple $(s,u,\chi)\in \MMM_t$ in terms of $G$-bundles
as follows. 

\ab First, identify ${\mathcal{P}}{\it{ic}}^\circ(\E)$,
the Picard variety of degree 0 line bundles on 
the elliptic curve $\E=\C^*\!/q^\Z$, with $\C^*\!/q^\Z$.
Let ${\cal L}_t\in {\mathcal{P}}{\it{ic}}^\circ(\E)$
denote the degree 0 line bundle  corresponding to 
the image of the complex number $t$ under the projection:
$\C^*\onto \C^*\!/q^\Z={\mathcal{P}}{\it{ic}}^\circ(\E)$.
Further,
given a principal semisimple  $G$-bundle $P\in\MM(\E, G)^{ss}$,
write $\g_{_P}$ for the associated
vector bundle corresponding to the adjoint representation in $\g=\Lie
G$. We call a pair $(P, x)$, where $x$ is 
 a regular section
of $\g_{_P}\otimes {\cal L}_t\,,$ a {\it
Higgs bundle}. Let $\Higgs(\E,G)^{^{\sf{nil}}}$ be the moduli
space of isomorphism classes
of triples $(P, x, \chi)$, where $P\in\MM(\E, G)^{ss}$, $x$ is 
a {\it nilpotent} regular section of 
$\g_{_P}\otimes {\cal L}_t\,,$ and  $\chi$ is an irreducible
admissible representation of
$Aut(P, x)/Aut^{\circ}(P,x)$,
the component group of the group of automorphisms of the Higgs bundle
$(P,x)$.

\ab
We claim that there is the following canonical
bijection, that should be thought of as
 a `$t$-deformation' of the bijection (5.14):
$$
\Big\{q{\mbox{\footnotesize{\it-conjugacy classes of triples }}}
(s,u,\chi)\in
\MMM_t\Big\}\quad\longleftrightarrow\quad\Higgs(\E,G)^{^{\sf{nil}}}
\eqno(6.4)
$$
The bijection assigns to a triple $(s,u,\chi)\in
\MMM_t\,$ the triple  $(P,x,\chi)$, where 
$P$ is the semi-simple bundle attached to the semisimple element $s$,
see (5.8.1), and $x$ is a  section
of $\g_{_P}\otimes {\cal L}_t\,$
arising from the loop:
$\,\log(u)\in \Lie(G_{q,s})\subset\g((z))\,,$ which is well-defined since
$u$ is unipotent. It is straightforward to see, using the
results of 
\S5, that this assignment sets up a bijection as in (6.4).
\bigskip

\ab
Based on the similarity with Theorem 5.16, we propose the following 
double-affine
version of the Deligne-Langlands-Lusztig conjecture
for affine Hecke algebras (proved in [KL2], see also [CG]).
\bigskip

{\bf Deligne-Langlands-Lusztig Conjecture for Cherednik algebras 6.5.} 

{\it If $q$ and $t$ are not roots of unity, then
there exists a canonical bijection between the set of
(isomorphism classes of) simple objects of the category
${\mathcal{M}}_{\!_{q,t}}(\hdd\,,\,\C\X)$
and the set of $\epq$-conjugacy classes in $\MMM_t$.}
\medskip

\ab Another evidence in favor
of Conjecture 6.5 comes from the  result of Garland-Grojnowski
 announced  [GG]. Garland-Grojnowski gave
a construction of the double affine Hecke algebra
in terms of equivariant $K$-theory of some infinite-dimensional
space. 
Modulo several  `infinite dimesionality' 
difficulties, Conjecture 6.5 might have been deduced from the
$K$-theoretic realization using the
nowadays standard techniques, see  [CG]. 
Unfortunately, the difficulties arising from `infinite dimesionality'
are extremely serious, and at the moment we do not see any way
to overcome them. That might look strange since
infinite
dimensional spaces almost never appear in [GG] explicitly, and the
authors of [GG] always avoid them by working
 with their finite-dimensional approximations.
This is deceptive, however, because in order to apply the
standard  techniques of [CG], one has to reformulate the constructions
of [GG] in manifestly infinite-dimensional terms
involving, in particular, equivariant $K$-theory with respect
to an {\it infinite-dimensional} group like $G((z))$
(as opposed to the $T$-equivariant $K$-theory used in [GG]).
Unfortunately, such a theory does not exist at the moment, for instance,
it is not even clear what should be the corresponding 
equivariant $K$-group of a point.
As a consequence, the crucial "localization at fixed points"
reduction does not apply. 
\bigskip

\section{Operator realization of the Cherednik algebra.}

\ab For each
$\mu \in  \Y=X_*(T)$, we introduce a $\mathbf{q}$-{\it shift
operator}, see \cite{kir},
$\,\dd^\mu: \C(T) \to \C(T)\,$ by letting it act 
by the formula
$$
(\dd^\mu f)(t)= f({\eps}^{2\mu}\cdot t)\quad,\quad
t\in
T.
$$ 
The operators of multiplication by $e^\la\,,\, \la\in \X=X^*(T),\,$ and
$\dd^\mu\,,\, \mu
\in \Y\,,$ satisfy
the commutation relation:
$\;\dd^\mu\ccirc e^\lambda =
{\eps}^{2\langle\lambda,\mu\rangle}\cdot e^\lambda\ccirc \dd^\mu\,.$
Thus, these operators generate an algebra isomorphic to
the quantum torus ${\hh=\hh(\X\oplus\Y).}$

\ab We consider $\C$-linear endomorphisms of $\C(T)$ of the form
$$
   h = \sum_{w\in W,\,\mu \in \Y}\; h_{w,\mu}\cdot  \dd^\mu  \cdot [w]\;
:\;\;\; f\mapsto \sum_{w\in W, \mu \in \Y}\; h_{w,\mu}\cdot 
\dd^\mu({}^{w\!}f)\,, \eqno(7.1)
$$
where $h_{w,\mu}\in\C(T)\,,$ and $f \mapsto {}^{w\!}f$
 denotes the natural action of $w\in W$ on $\C(T)$.
The set of all such operators forms an associative algebra
$\hh_{_{{\sf frac}}}[W]$, isomorphic to a smash product of 
the group algebra $\C[W]$ with the algebra $\hh_{_{{\sf frac}}}$
of difference
operators on $T$ with rational coefficients.
Observe that the algebra
$\hh_{_{{\sf frac}}}$  is a slight enlargement of the quantum torus
algebra $\hh=\hh(\X \oplus \Y)$.
The difference between $\hh_{_{{\sf frac}}}$ and $\hh$ is that
we are allowing coefficients of difference operators to be rational,
not just polynomial.

\ab In
\cite{GKV}, the  Cherednik algebra $\hdd$ has been realized as
a subalgebra of the smash product algebra $\CqepsT  
\st\, \C [\Waff  ].$ 
Theorem 7.2 below gives
a similar description of $\hdd$ with the affine Weyl group $\Waff$ being
replaced by the finite group $W$ at the cost of replacing
the commutative algebra $\CqepsT$ by a non-commutative algebra
of finite-difference operators. Note that the smash product
algebra $\CqepsT  \st \C[\Waff  ]$ may
be viewed  as a smash product $\hh_{_{{\sf frac}}}[W]=
\hh_{_{{\sf frac}}}  \st \C [W]$.

\ab For a root $\beta\in \Delta$ let $\epsilon(\beta)=1$ if
$\beta \notin R^+ ,$ $\epsilon(\beta)=0$ otherwise.
Given $\tau\in \C,$ and a root $\alpha \in \Delta,$ let 
$T_{\alpha,\tau}$ be the divisor in $T$ given by the condition 
$e^{\alpha}=\tau.$ Let $T_\alpha=T_{\alpha,1}.$
The main result of this section is the following
operator description of the double-affine  Hecke algebra $\hdd$,
similar to the description of the affine Hecke algebra given in [GKV].
\medskip

{\bf Theorem 7.2.} {\it The algebra $\hdd$ is isomorphic to the subalgebra
of $\hh_{_{{\sf frac}}}[W]$ formed by all the elements
$h = \sum_{w\in W, \mu \in \Y}\; h_{w,\mu} \dd^\mu  \cdot [w]$
whose coefficients $h_{w,\mu}\in\C(T)$ satisfy the following
conditions:}

{\bf{(7.2.1)}}$\;\;$  $h_{w,\mu}$ {\it is regular, except at divisors
$T_{\alpha,{\eps}^{2k}}\,,\, k\in \Z$, where they may have first order poles.}

{\bf{(7.2.2)}}$\;\;$ $ {\Res}_{_{T_{\alpha,{\eps}^{-2k}}}}(h_{w,\mu})
+{\Res}_{_{T_{\alpha,{\eps}^{-2k}}}}(h_{s_\alpha w,k\alpha +s_{\alpha}\mu})
=0,\quad\forall\alpha\in \Delta\,;
$

{\bf{(7.2.3)}}$\;\;$ {\it For each $\alpha\in \Delta^+$ the function  $h_{w,\mu}$
vanishes at the divisor $T_{\alpha, p}$ for the following values of $p$:}

$\hspace{16mm}$ {\small \parbox[t]{125mm}{
$\dis p={\eps}^{2k}\zht^{-2}
\quad \mbox{if}\quad \langle\alpha, \mu\rangle < 0,\quad \mbox{and}\quad 
0 \le k \le |\langle\alpha , \mu\rangle + 1
-\epsilon(w^{-1}(\alpha))|$\newline
$\dis p=\zht^{-2}
\quad\,\,\,\quad \mbox{if}\quad \langle\alpha, \mu\rangle  = 0,
\quad \mbox{and}\quad \epsilon(w^{-1}(\alpha))=1$\newline
$\dis p={\eps}^{-2k}\zht^2
\quad \mbox{if}\quad \langle\alpha, \mu\rangle  > 0,
\quad \mbox{and}\quad 
 1 \le k \le |\langle\alpha , \mu\rangle  
-\epsilon(w^{-1}(\alpha))|$. 
}}
\bigskip

\ab The rest of this section is devoted to the proof of this Theorem,
which will be based on the `zero-residue' construction of the algebra
$\hdd$ given in [GKV].

\ab  Recall that the affine root system has the subset of real roots
$\Raffreal$ which are the affine roots whose restriction to $T\subset
\dT$ is nonconstant. 
Given $\tau\in \C,$ and a root $\gamma
\in \Raf$,  let 
$\dT_{\gamma,\tau}$ be the divisor in $\dT$ given by the condition 
$e^{\gamma}=\tau.$ Let ${\dT}_\gamma={\dT}_{\gamma,1}.$
According to [GKV], each operator $f\in \hdd$ is written as
$$
  f = \sum\nolimits_{w\in \Waff}\; f_w [w] \eqno(7.3)
$$
and the  coefficients $f_w$  satisfy certain zero-residue  conditions. 
The conditions in \cite{GKV} on the coefficients $f_w$ for
$w\in \Waff$ are as follows:

\begin{quote}
{\bf{(7.4.1)}}$\;\;$ Each $f_w$ has at worst first order poles at the divisors
$\TKMgamma,$ for $\gamma$ a real root, and is otherwise regular.\par\noindent
{\bf{(7.4.2)}}$\;\;$ For each $w\in \Waff$ and real root $\gamma,$ we have:
\[{\Res}_{\TKMgamma}(f_w)+ {\Res}_{\TKMgamma}(f_{s_\gamma w})=0.\]
{\bf{(7.4.3)}}$\;\;$ The function $f_w$ vanishes on $\TKMaqm$ whenever $\gamma$ is
a positive real root and $w^{-1}(\gamma)<0.$
\end{quote}
\medskip

Since $\Waff = \Y\rtimes W$  we can rewrite the  expression on the RHS
of (7.3)
as
\[
  f = \sum_{w\in W,\, \mu \in \Y } h_{w,\mu}\;\; [\mu \cdot w]\,.\]
The $h_{w,\mu}$ satisfy certain zero-residue conditions arising from
(7.4.1)--(7.4.3).
We are now  going to translate each of the conditions $7.4.(i)\,,\, i=1,2,3,$
into the corresponding conditions of
Theorem $7.2.(i)$.\medskip

\ab ${\mathbf {(7.4.1)\enspace\Longrightarrow\enspace(7.2.1)}}
\quad$ The real roots are of the form $\alpha + k\delta$, $\alpha \in
\Delta$, $k\in \Z$. We take $e^\delta = {\eps}^{2}$. Then the condition
$e^{\alpha + k\delta}=1$ defines the divisor $T_{\alpha,{\eps}^{-2k}}$.
\medskip

\ab ${\mathbf {(7.4.2)\enspace\Longrightarrow\enspace(7.2.2)}}
\quad$ 
In $\Waff$, we have the formula:
$\;\dis s_{\alpha + k\delta} {\eps}^\nu w = 
{\eps}^{\alpha + s_\alpha \mu} s_\alpha w\,,$
where $\alpha + k\delta$ is a real root, $s_{\alpha + k\delta}$ is
the corresponding reflection, $w\in W$, and for $\nu \in Y$, ${\eps}^\nu$
represents the corresponding translation in $\Waff$. Then ${\mathbf{ (7.4.2)
}}$, in the case $\beta=\alpha + k\delta$, gives ${\mathbf{(7.2.2)}}$.
\medskip

\ab ${\mathbf {(7.4.3)\enspace\Longrightarrow\enspace(7.2.3)}}
\quad$ We know that $h_{w,\mu}$ vanishes at the divisor $\TKMgqm$
whenever $\gamma \in \Raffrepos$ and $(\mu w)^{-1}(\gamma)\not\in 
\Raffrepos$. For $\alpha\in\Delta^+$, 
the divisor $T_{\alpha, p}$ can arise from three
types of roots in $\Raffrepos$: (i) $\;\alpha + k\delta\,,\, k > 0\,,\;$
(ii) $\;\alpha\,,\;$ (iii) $\;-\alpha + k\delta\,,\, k > 0$. In case (i), we get:

$$(\mu w)^{-1}(\alpha + k\delta)=w^{-1}(\alpha) + (k+\langle\alpha,\mu\rangle)\delta\,.
\eqno(7.6)$$
Hence, $h_{w,\mu}$ vanishes at ${\TKM}_{\alpha + l\delta,{\zht}^{-2}}$, provided
$\langle\alpha,\mu\rangle < 0$ and $l=0,1,\ldots, -\langle\alpha,\mu\rangle$ if
 $\;\epsilon(w^{-1}\alpha)=1,\;$ and for $l=0,1,\dots, -\langle\alpha,\mu\rangle - 1$
if $\;\epsilon(w^{-1}\alpha)=0$. The equation for the divisor:
${\TKM}_{\alpha + l\delta,{\zht}^{-2}}= T_{\alpha,{\eps}^{2l}{\zht}^{-2}}$,
yields the first case. Cases (ii) and (iii) follow similarly,
using (7.6).\sq\bigskip

\ab We would like to propose a characterisation of the Cherednik algebra
similar to
the characterisation of  the affine Hecke algebra given 
in [GKV, Theorem 2.2].
To this end, let $M_{_{{\sf frac}}}$ denote a rank one vector space over $\C(\dT)$ with
generator $\m$. For each $i=0,1,\ldots,l,$ define an action of 
$s_i\in \Waff$ on $M_{_{{\sf frac}}}$ by the formula:
$$
\hat{s}_i: \; f\cdot\m \;\mapsto \;s_i(f)\cdot \frac{\eps^{-1}\cdot e^{\alpha_i/2}-
\eps\cdot  e^{-\alpha_i/2}}{\eps^{-1}\cdot e^{-\alpha_i/2}-
\eps\cdot  e^{\alpha_i/2}}\cdot\m\quad,\quad f\in \C(\dT)
\eqno(7.7)
$$
It is easy to see that the assignment: $s_i\mapsto \hat{s}_i\,,\,
i=0,1,\ldots,l,$ extends to a representation of
the affine Weyl group $\Waff$ on $M_{_{{\sf frac}}}$. This way one makes
$M_{_{{\sf frac}}}$ a module over the smash-product algebra $\C[\Waff]\st\C(\dT)$.
Let $M\subset M_{_{{\sf frac}}}$ be the
free $\C[\dT]$-submodule generated by $\m$.
It is straightforward to verify that $M$ is stable under the
action of the elements:
$\hT_i\in \C[\Waff]\st\C(\dT)\,,\,
i=0,1,\ldots,l,$  defined by formula (6.1).
\medskip

{\bf {Question.}} How to reformulate  (7.7) in terms of difference
operators ?
\bigskip 

\section{Spherical subalgebra}

\ab Set $W(\zht) := \sum_w\, \zht^{2l(w)}\,,$ and let
$\,\ee=\frac{1}{W(\zht)}\cdot \sum_w\, \zht^{l(w)}\cdot \T_w\,\in \HH$
be the central idempotent, see [KL], corresponding  to   the   1-dimensional
$\HH$-module:  $\T_w\mapsto  \zht^{l(w)}$. 
We identify $\HH$ with a subalgebra of $\hdd$ in a natural way, and
regard $\ee$ as an element of $\hdd$.
 Write    $\ehe$   for    the    subalgebra
$\ee\cdot\hdd\cdot\ee\sset\hdd$, which we call the {\it spherical
subalgebra}. Unlike the case of affine Hecke algebra, the subalgebra
  $\ehe$  is not commutative.

\ab   We now give an explicit `operator' description of the Spherical
subalgebra. Let $s_1,\cdots, s_l$ be the simple reflections in $W$.

\medskip

{\bf Theorem 8.1.} {\it An element $h=\sum_{\mu\in \Y\,,\,w\in W}\,
h_{w,\mu}\cdot
\dd^\mu\,[w]
\in\hdd$, see Theorem 7.2, belongs to $\ehe$ if and only if,
for any $i=1,\ldots,l,$ we have:}

\ab  \vi $\hphantom{x}\qquad\dis h_{s_iw,s_i\mu}=s_i(h_{w,\mu})\,;$

\ab \vii$\hphantom{x}\qquad\dis h_{ws_i,\mu}=h_{w,\mu}\cdot
 \dd^\mu w({{\zht^2e^{\alpha_i}-1}
\over {e^{\alpha_i}-\zht^2}}) \,.$
\medskip

{\sl  Proof.} We first consider the $SL_2$ case. Since $\ee$ is idempotent,
$\ehe$ is the $1$--eigenspace for the left and right actions of $\ee$.
From the formula for the operator $T_\alpha$ given by (6.1),
we can write $\ee= a\cdot s_\alpha + b\cdot 1$ with 
$$
 a = {{\zht^2 e^{\alpha} - 1}\over {(1+{\zht}^2)(e^{\alpha}-1)}}\quad;
\quad b = {{ e^{\alpha} - \zht^2}\over {(1+{\zht}^2)(e^{\alpha}-1)}}\,.
$$
It is straightforward to check that: $s_\alpha (a)=b$ and $a= 1 - b$.
Write: $h = \sum_{n\in \Z}\, (h_n \dd^{n\rhocheck}s_\alpha
+ g_n \dd^{n\rhocheck}).\,$
Then, the equation $\ee\cdot h=h$ implies: $\, h_n=s_\alpha (g_{-n})\,.$
Similarly, the equation $h= h\cdot \ee$ implies: 
$$
 h_n=g_n\cdot ({{\zht^2{\eps}^{2n}e^{\alpha}-1}\over
{{\eps}^{2n}e^{\alpha}-\zht^2}})\;.
$$

\ab The general case follows from the fact that if $M$ is a left or right
 module for the finite Hecke algebra, then $m\in \ee\cdot M$ if and
only if $T_i\cdot m=\zht\cdot  m$ for all $i=1, \dots, l\,.\quad\square$
\medskip

{\bf Remark.}\ab We note that
$\;
({{\zht^2{\eps}^{2n}e^{\alpha}-1}\over
{{\eps}^{2n}e^{\alpha}-\zht^2}})
 = \dd^{n\rhocheck}({{\zht^2 e^{\alpha}-1}\over {e^{\alpha}-\zht^2}}).
$
Moreover, the expression $\;({{\zht^2 e^{\alpha}-1}\over {e^{\alpha}-\zht^2}})\,$
appears in work of Drinfeld on affine quantum groups [Dr], and has been
interpreted as a characteristic class in [GV].
\bigskip

\ab    We can deform
the idempotent $\ee$ to obtain a  family of
 idempotents $\ee_v \in \hdd_v$ as follows.
Given $w\in W$, write its reduced  decomposition:
$w=s_{i_1}\cdots s_{i_l}$.  Put: 
$\,\T_{w,v} :=\T_{i_1,v}\cdot \ldots\cdot \T_{i_l,v}\in \hdd_v\,.$ It is
standard to show that this element
is independent of the choice of reduced decomposition of $w$.
 We let: ${\mathbf{y}}:=\zht -v(\zht -
\zht^{-1})$,
and put:
$$
W(\zht,v):=\sum\nolimits_{w\in W}\;   \left( {{\zht}\over {{\mathbf{y}}}}\right)^{l(w)}
\quad,\quad
\ee_v :={1\over {W(\zht,v)}}\cdot\sum\nolimits_{w\in W}\; 
\frac{1}{{\mathbf{y}}^{l(w)}}\T_{w,v}
$$
It is easy to check that $\T_{i,v}\cdot \ee_v = \zht\cdot  \ee_v,$ and
to derive from this that $\ee_v$ is an idempotent.
Using the family of idempotents: $\ee_v\in \hdd_v$ we define
a family of spherical subalgebras: $\,\ehev\subset \hdd_v.\,$
Note that for $v=0$ we have: $\,\ehez\cong  \hw.\,$
  
\ab We  would like   to  deform   the  representations of $\hw$
constructed in \S4.   Assume from now on that both
$q$ and $t$ are {\sf not roots of unity}.\medskip

{\bf  Deformation Conjecture 8.2.}  {\it If $\Delta$ is the root system
of type ${\mathbf{A_n}}$ then,
for  any $\lambda \in \La\,,\,
\chi  \in {\widehat W}^\la$,   the  simple  
$\hh[W]$-module $Z_\chi$,  see
  Proposition 3.5, can be  deformed, for each $v\neq 0$,
 to a simple object of the category
  ${\mathcal{M}}_{\!_{q,t}}(\hdd_{\!v}\,,\, \C\X)$.}  \medskip

\ab 
The most trivial representation of $\hw$ is the one
in the space of invariant polynomials:
 ${\C  [T]}^W=\ee_{_0} \C  [T]\,$. This space is easily deformed to $\ee_v \C [T]$,
yielding a representation of $\ehev$.

\ab Next, we deform the sign representation.
To this end, define the {\it Iwahori-Matsumoto involution} $\Xi$
on the algebra
 $\hddv$ by
the formulas:
$$\Xi(\T_{i,v})=v(\zht-\zht^{-1}) - \T_{i,v}, 
\qquad \Xi(e^\mu)= e^{-\mu}\,.$$
Let $\,\alt_v=\sum_{w\in W}\;(-1)^{\ell(w)}\cdot\T_{w,v} \in \HH$
denote the standard anti-symmetriser.
It is easy to check that $\Xi$ is an algebra automorphism
such that: $\Xi(\ee_v)= \alt_v$.
Hence, the Iwahori-Matsumoto involution gives an algebra isomorphism
$\,\Xi: \ehev \iso \alt_v\hdd_v\alt_v$. Composing this
 isomorphism with the natural $\alt_v\hdd_v\alt_v$-action
on the space $\alt_v\C  [T]$ we get a new `sign-representation'
of $\ehev$. 
For $v=0$, for instance,
this gives the representation of $\hw$ on the space: $\,\C [T]^{sgn}
:=$
$ \{P\in \C [T]\;\;\Big|\;\; w(P)=sgn(w)\cdot P\;,\; \forall w\in W \}\,.$
Thus, we have constructed a family of representations of $\ehev$,
a deformation of the representation of $\hw$ corresponding to the
sign representation.
It is likely that, for $\g={\mathfrak{g}\mathfrak{l}}_n$, 
an extention of this construction
to more general Young symmetrisers, would allow to deform 
simple $\hw$-modules corresponding to other representations
of the Symmetric group, cf. Theorem 4.3.\medskip

\ab 
Provided Conjecture 8.2 holds we expect, moreover, that
(in the non root of unity case) all simple objects of
${\mathcal{M}}_{\!_{q,t}}(\hdd_{\!v}\,,\,\C\X)$ can be
obtained by deformation of simple $\hh[W]$-modules $Z_\chi$
(for type ${\mathbf{A_n}}$).
 \bigskip

{\bf Remark 8.3.} There are natural functors, cf. Proposition 4.1:
\[
{\bold F}: \hdd\mod  \sto   \ehe\mod \quad,\quad M 
\mapsto \ee\hdd\;\bigotimes\nolimits_{_{\hdd}}\;M\]
\[{\bold I}: \ehe\mod  \sto\hdd\mod \quad,\quad 
N \mapsto
\hdd\ee\;\bigotimes\nolimits_{_{\ehe}}\;N \]
However, unlike the situation considered in \S4,
these functors do not give rise to Morita equivalence, in general,
because the natural maps:
$$\hdd\ee\;\bigotimes\nolimits_{_{\ehe}}\;\ee\hdd \too\hdd\quad,\quad
\ee\hdd\;\bigotimes\nolimits_{_{\hdd}}\;\hdd\ee \too \ehe$$
generally fail to be injective. We expect that the functors
${\bold F}$ and ${\bold I}$ do provide a Morita equivalence
if $\zht$ and ${{\bold q}}$ are specialized to complex numbers
$t$ and $q$, respectively, such that $t^m
\neq q^n$, for any $(m,n)\in \Z^2\smallsetminus(0,0)$.
\bigskip

\section*{Springer correspondence for disconnected groups.}

In  this  Appendix  we  show  how to  extend  the  classical  Springer
Correspondence  to the  case  of not  necessarily connected  reductive
groups. First we recall briefly  (see [CG, Chapter 3] for details) the
situation for a general  \textit{connected} reductive group, such as the
group $H^{\circ}$ of \S5.

\ab  Let  ${\mathcal{N}} \subset  H^{\circ}$  be  the  subset  of  unipotent
elements  and  $\mathcal{B}$   the  variety  of  Borel subgroups  in
$H^{\circ}$.    The   subvariety    $\widetilde{\mathcal{N}}   \subset
\mathcal{B} \times \mathcal{N}$ of all pairs  $\{ (B_H, u) \enspace|\enspace u \in B_H
\}$,  provides  an  $H^{\circ}$-equivariant  smooth  resolution  $\pi:
\widetilde{\mathcal{N}} \to \mathcal{N}$, called the \textit{Springer
  resolution}. 

\ab  Denote   by    ${\mathcal{Z}}$   the   fiber    product   $\widetilde{\mathcal{N}}
\times_{\mathcal{N}}\widetilde{\mathcal{N}}$,   which   can  also   be
indentified with (cf. [CG])  the subvariety in $T^*(\mathcal{B} \times
\mathcal{B})$  given   by  the  union  of  conormal   bundles  to  the
$H^{\circ}$-orbits on  $\mathcal{B} \times \mathcal{B}$  (with respect
to the diagonal action). The  top Borel-Moore homology group $H({\mathcal{Z}})$ is
endowed with a structure of an associative algebra via the convolution
product (see  [CG]). Moreover, the set ${\mathbb{W}}  \subset H({\mathcal{Z}})$ of
fundamental classes  of irreducible components  of ${\mathcal{Z}}$, forms  a group
with respect  to the convolution product,  called the \textit{abstract
  Weyl group}, and $H({\mathcal{Z}})$ can  be identified with the group algebra of
$\mathbb{W}$. A particular choice of  a Borel subgroup $B_H \supset T$
identifies the  usual Weyl group  $W^{\circ}= N_{H^{\circ}}(T)/T$ with
$\mathbb{W}$ by  sending the class  of $n_w \in N_{H^{\circ}}$  to the
fundamental class  of the conormal bundle to  the $H^{\circ}$-orbit of
$(B_H, n_w B_H n_w^{-1}) \in \mathcal{B} \times \mathcal{B}$.

\ \  \ \ab  Consider a  unipotent orbit $\mathcal{O} \subset  \mathcal{N}$. 
The  top   Borel-Moore  homology  groups   of  the  fibers   of  $\pi:
\widetilde{\mathcal{N}} \to  \mathcal{N}$ over $\mathcal{O}$,  form an
irreducible local  system $L_{\mathcal{O}}$ on  $\mathcal{O}$ which is
equivariant  with  respect  to  ${\mathbb{W}} \times  H^{\circ}$  (the
action of  $\mathbb{W}$ in the fibers of  $L_{\mathcal{O}}$ comes from
the convolution  construction, cf [CG], and the  action of $H^{\circ}$
from    the    $H^{\circ}$-equivariance    of    $\pi$).     Decompose
$L_{\mathcal{O}}$  into  a  direct  sum of  irreducible  ${\mathbb{W}}
\times H^{\circ}$-equivariant  local systems $L_1,  \ldots, L_k$.  For
any representation  $\phi$ of $\mathbb{W}$  we can consider  the local
system  $I_{i}$ formed  by the  $\mathbb{W}$-invariants of  the tensor
product  $\phi^\vee  \otimes  L_i$.    It  turns  our  that,  for  any
irreducible  representation  $\phi$,   there  exists  a  unique  orbit
$\mathcal{O_{\phi}}$ and a unique  $L_{\phi} \in \{L_1, \ldots, L_k\}$
for  which the  local system  $I_{\phi}$, constructed  from  $\phi$ as
above,  is  non-zero.  Moreover,  such  $I_{\phi}$  is an  irreducible
$H^{\circ}$-equivariant local  system associated to  a ``admissible
representation''  (in the sense  of Definition  5.2) of  the component
group of the centralizer $Z_u$ of a point $u \in {\mathcal{O}}$. Below
we  will use  the  language  of equivariant  local  systems (which  is
equivalent to the language of admissible representations).

\ab  The      \textit{Springer      correspondence}      $\phi      \mapsto
({\mathcal{O}}_{\phi}, I_{\phi})$ gives a bijection between the set of
irreducible representations of its Weyl group $\mathbb{W}$ and the set
of pairs  $({\mathcal{O}}, I)$ where $\mathcal{O}$ is  a unipotent orbit
of   $H^{\circ}$  and   $I$  is   a   certain  $H^{\circ}$-equivariant
irreducible local system on  $\mathcal{O}$ coming from a admissible
representation of $Z_u/Z_u^{\circ}$.

\medskip \ab  We proceed to representation theory
of   the    `Weyl   group'   $W_H  =    N_H(T)/T$   of   a
\textit{disconnected}  reductive  group  $H$  (see Prop. 5.13).

\bigskip
\textbf{Lemma A.1} \emph{A choice of a  Borel subgroup $B_H  \supset
  T$ in $H^{\circ}$ identifies the Weyl group $W_H$ with the
  semidirect product $W^{\circ} \rtimes (W_H/W^{\circ})$. Moreover, one
  has a canonical isomorphism
  $H/H^{\circ}  \simeq W_H/W_H^{\circ}$.}

\noindent
\textit{Proof.} Consider  the subgroup $N'(T):=  N_H(B_H) \cap
N_H(T)$. Then the embedding $N'(T) \subset H$ induces the isomorphisms
$$
H/H^{\circ} \simeq N'(T)/T \simeq W_H /W^{\circ} 
$$
Since $N'(T)$ is a subgroup of $N_H(T)$, we obtain an embedding
$W_H/W^{\circ} \subset W_H$. Now the assertion of the Lemma follows.
\hfill $\square$

\bigskip
\noindent
\textbf{Remark.} A different choice Borel subgroup $B'_H$ containing $T$
gives a conjugate embedding $w(W_H/W^{\circ}) w^{-1} \subset W_H$,
where $w \in W^{\circ}$ is the unique element which conjugates $B_H$
into $B'_H$.

\ab  Note that $H$  acts on $\mathcal{N}$ and
on $\widetilde{\mathcal{N}}$. 
In  particular $H$ permutes  the irreducible  components of
${\mathcal{Z}}$.
That
induces  an $H/H^{\circ}$-action  on  $\mathbb{W}$   by  group
automorphisms.

\bigskip
\noindent
\textbf{Proposition  A.2.}    \emph{The  isomorphism  ${\mathbb{W}}  =
  W^{\circ}$  (depending on  the choice  of $B_H$)  and  the canonical
  isomorphism $H/H^{\circ}  \simeq W_H/W^{\circ}$, identify  the above
  action  of  $H/H^{\circ}$ on  ${\mathbb  W}$,  with the  conjugation
  action of $W_H/W^{\circ}$ on $W^{\circ}$ arising from Lemma A.2.}

\medskip
\noindent
\textit{Proof.}  It  suffices to replace the pair  $(H, H^{\circ})$ by
$(N'(T),  T)$. Let $n_w$  be a  lift to  $N_{H^{\circ}}$ of  a certain
element $w \in W$, and let  ${\mathcal{Z}}_w$ be the cotangent bundle to the orbit
of  $(B_H, n_w  B_H  n_w^{-1}) \in  \mathcal{B}  \times \mathcal{B}$.  
Similarly, let $n_{\sigma} \in N'(T)$  be a lift of an element $\sigma
\in  W_H  /W^{\circ}$.  Denote  $\sigma  w  \sigma^{-1} \in  W^{\circ}
\subset W_H$  by $w^{\sigma}$,  then $n_{w^{\sigma}} =  n_{\sigma} n_w
n_{\sigma}^{-1}$ is a lift of $w^{\sigma}$ to $N_{H^{\circ}}(T)$.
 
\ab    By definition of $N'(T)$  the element $n_{\sigma}$ normalizes $B_H$. 
  Hence $n_w$ sends $(B_H, n_w B n_w^{-1})$ to $(B_H, n_{w^{\sigma}} B_H
  n_{w^{\sigma}}^{-1})$.    Thus,    $n_{\sigma}   \cdot   {\mathcal{Z}}_w   =
  {\mathcal{Z}}_{w^{\sigma}}$.$\quad\square$

\bigskip

\ab  Now we recall the basic facts of Clifford
theory (cf. [Hu]) which apply to any finite group $W_H$ and its
normal subroup $W^{\circ}$, not necessarily arising as Weyl groups.

\ab  The group $W_H$ acts by conjugation on the set $\widehat{W}^{\circ}$ of
irreducible representations of $W^{\circ}$. Let ${\mathcal{V}}_1 \dots
{\mathcal{V}}_k$ be the orbits of its action. For any irreducible
representation $\psi \in \widehat{W}_H$ we can find an orbit
${\mathcal{V}}_{i(\psi)}$ and a positive integer $e$, such that the
restriction of $\psi$ to $W^{\circ}$ is isomorphic to a multiple of the
orbit sum: 
$$
\psi|_{W_0} \simeq e \cdot \left(\sum\nolimits_{\phi  \in {\mathcal{V}}_{i(\psi)}
  } \phi \right)
$$
Fix $\phi \in {\mathcal{V}}_{i(\psi)}$ and consider the subset
$(\widehat{W}_H)_{\phi} \subset \widehat{W}_H$ of all representations whose
restriction to $W^{\circ}$ contains an isotypical component isomorphic to
$\phi$ (and hence automatically all representations in the orbit of
$\phi$). Obviously, $\widehat{W}_H$ is a disjoint union of
$(\widehat{W}_H)_{\phi_i}$, where $\phi_i \in {\mathcal{V}}_i$ is any
representative of the orbit ${\mathcal{V}}_i$.

\ab  To   study  $(\widehat{W}_H)_{\phi}$   we  consider   the  stabilizer
$W^{\phi} \subset  W_a$ of  $\phi \in \widehat{W}^{\circ}$.   Then by
Clifford  theory (cf. [Hu]),  the induction  from $W^{\phi}$  to $W_H$
establishes  a  bijection  between  $(\widehat{W}^{\phi})_{\phi}$  and
$(\widehat{W}_H)_{\phi}$.  Moreover, any linear representation $\chi \in
(\widehat{W}^{\phi})_{\phi}$ is isomorphic  to the tensor product $p_1
\otimes  p_2$  of two  \textit{projective}  representations $p_1$  and
$p_2$ (cf.  [Hu]) such that

\ab  \vi $p_1(x) = \phi(x)$, $p_2(x) = 1$ if $x \in W^{\circ}$

\ab  \vii $p_1(gx) = p_1(g) \phi(x)$ and $p_1(xg) = \phi(x) p_1(g)$ if $ x
\in W^{\circ}, g \in W_H$.

\smallskip

\ab  Thus, $p_2$  is a projective representation of  $W_H/ W^{\circ}$ which
plays the role of the multiplicity  space of dimension $e$ in terms of
the above formula for $\psi|_{W^{\circ}}$. The second  condition means
that the projective cocycle of $p_1$ is in fact lifted from $W^{\phi} /
W^{\circ}$. From now on we will fix the decomposition $\chi = p_1
\otimes p_2$. 

\ab Let $\widehat{W}_H$ denote the set of isomorphism
classes of irreducible representations of $W_H$. Further, for any 
unipotent conjugacy class ${\mathcal{O}} \subset H$, let
${\mathsf{Admiss}}({\mathcal{O}})$ stand for the set of (isomorphism
classes of) irreducible admissible (in the sense
of Definition 5.2) $H$-equivariant local systems on ${\mathcal{O}}$.

\bigskip
\noindent
\textbf{Theorem A.3.} \emph{There exists a bijection between the
following sets}:
$$\widehat{W}_H\enspace\longleftrightarrow\quad
\Big\{
{{\mbox{\footnotesize{unipotent conjugacy
class}}}\atop{\mbox{\footnotesize
{${\mathcal{O}} \subset H$ and $\alpha\in {\mathsf{Admiss}}({\mathcal{O}})$
}}}}\Big\}
$$
\medskip
\noindent
\textit{Proof.} Take an irreducible representation $\rho$ of $W_H$
and let $\phi \in \widehat{W}^{\circ}$ be an irreducible
subrepresentation of $\rho|_{W^{\circ}}$. By Clifford theory $\rho$ is
induced from a certain representation $\chi \in
(\widehat{W}^{\phi})_{\phi}$ as above.

\ab  Recall that by Springer Correspondence for $W^{\circ}$ the irreducible
representation  $\phi$ gives  rise  to a  unipotent $H^{\circ}$  orbit
$\mathcal{O}_{\phi}$  together with  an  $H^{\circ}$-equivariant local
system $I_{\phi}$. We will show how to construct from $\chi = p_1
\otimes p_2$ the corresponding local system for $H$.

\ab  As a  first  step,  we   will  construct  a  certian  local  system
$\tilde{I}_{\chi}$  on  $\mathcal{O}_{\phi}$.   This local  system  is
equivariant with  respect to the  subgroup $H^{\phi} \subset  H$ which
corresponds to $W^{\phi} \subset W_H$ via the isomorphism $H/H^{\circ}
= W_H /W^{\circ}$ of Lemma A.1  (it is easy to prove using Proposition
A.2,  that $H^{\phi}$ is  the subgroup  of all  elements in  $H$ which
preserve  the orbit $\mathcal{O}$  and the  local system  $I_{\phi}$). 
Then,  imitating the  induction  map $(\widehat{W}^{\phi})_{\phi}  \to
(\widehat{W}_H)_{\phi}  $  we  will  obtain an  $H$-equivariant  local
system   on   the    unique   unipotent   $H$-orbit   which   contains
$\mathcal{O}_{\phi}$ as its connected component.

\ab  We fix the choice of  Borel subgroup $B_H$ and, in particular, the
isomorphisms: ${\mathbb{W}} \simeq W^{\circ}$
and  $W_H \simeq W^{\circ}
\rtimes (W_H/W^{\circ})$.

\bigskip
\noindent
{\bf {Step 1.}}\quad
Recall that the  Springer resolution $\pi: \widetilde{\mathcal{N}} \to
\mathcal{N}$ is $H$-equivariant. It  follows from the definitions that
there exists an action of $H^{\phi}$ on the total space of $L_{\phi}$,
which  extends  the  natural  action  of  $H^{\circ}$.   However,  the
extended action does not commute with the $W^{\circ}$-action any more.
Instead, it  satisfies the identity $h(w  \cdot s) =  h(w) \cdot h(s)$
where $h \in  H^{\phi}$, $w \in W^{\circ}$ and $s$  is a local section
of  $L_{\phi}$.   We  would  like  to  use  the  formula  $I_{\phi}  =
(\phi^\vee    \otimes    L_{\phi})^{W^{\circ}}$    to    define    the
$H^{\phi}$-equivariant structure on $I_{\phi}$.   To that end, we have
to construct an action of  $H^{\phi}$ on $\phi^\vee$ which agrees with
the $W^{\circ}$-action  in the  same way as  before.  This is  done by
using   the   composition   $H^{\phi}   \to  H^{\phi}/   H^{\circ}   =
W^{\phi}/W^{\circ} \hookrightarrow  W_H$ and the  projective action of
$W_H$  on   $p_1^\vee$  coming  from  Clifford   theory  (recall  that
$p_1^\vee$  extends  the $W^{\circ}$-action  on  the  vector space  of
$\phi^\vee$).   By Proposition  A.2  the two  actions  of $H^\phi$  on
$W^{\circ}$  coincide, hence  the projective  action of  $H^{\phi}$ on
$p_1^\vee  \otimes L_{\phi}$  indeed satisfies  $h(w \cdot  s)  = h(w)
\cdot  h(s)$. Consequently,  the  local system  $I_{\phi} =  (p_1^\vee
\otimes   L_{\phi})^{W^{\circ}}$  carries   a  projective   action  of
$H^{\phi}$.

\ab  Now we tensor the local system $I_{\phi}$ with the vector space of the
projective representation  $p_2^{\vee}$. Since $H^{\phi}$  acts on the
vector space  of $p_2^{\vee}$ via  the same composition  $H^{\phi} \to
H^{\phi}/H^{\circ} \simeq W^{\phi}/W^{\circ} \hookrightarrow W_H$, the
tensor    product   $p_2^{\vee}    \otimes   I_{\phi}$    carries   an
\textit{apriori} projective action  of $H^{\phi}$.  However, since the
projective cocyles  of $p_1$ and  $p_2$, well-defined as  functions on
$W^{\phi}  /W^{\circ}$, are  mutually inverse,  the same  can  be said
about the projective cocycles  of the $H^{\phi}$-actions on $I_{\phi}$
and  $p_2^{\vee}$.   Therefore,  these  cocycle cancel  out  giving  a
\textit{linear}  action  on  the  tensor  product.   This  means  that
$p_2^{\vee}   \otimes   I_{\phi}$  is   given   the   structure  of   an
$H^{\phi}$-equivariant    local    system,    to   be    denoted    by
$\tilde{I}_{\chi}$.

\medskip
\noindent
{\bf{Step 2.}}\quad
Next  we   consider  a  larger   subgroup  $\widehat{H}^{\phi}$  which
preserves   the   unipotent   orbit  $\mathcal{O_{\phi}}$,   but   not
necessarily the local system $I_{\phi}$.  It is easy to check that the
composition $\widehat{H}^{\phi} \times_{H^{\phi}} \tilde{I}_{\phi} \to
\widehat{H}^{\phi}     \times_{H^{\phi}}     \mathcal{O}_{\phi}    \to
\mathcal{O}_{\phi}$ defines  an $\widehat{H}^{\phi}$-equivariant local
system    $\hat{I}_{\chi}$   over    $\mathcal{O}_{\phi}$.    Finally,
$I_{\chi}=  H  \times_{\widehat{H}^{\phi}}  \hat{I}_{\chi}$ is  an
$H$-equivariant  local system over  $H \times_{\widehat{H}^{\phi}}
\mathcal{O}_{\phi}$.  Note  that the latter  space is nothing  but the
union  of those  unipotent $H^{\circ}$-orbits  which are  conjugate to
each  other with  respect  the the  larger  group $H$,  i.e. a  single
unipotent $H$-orbit.
  It is easy  to check that the assignment:  $\chi \mapsto I_{\chi}$,
$\chi \in (\widehat{W}^{\phi})_{\phi}$ together with the decomposition:
$\widehat{W}_H     =      \bigcup_{\phi_i     \in     {\mathcal{V}}_i}
(\widehat{W}_H)_{\phi}$     and     the     induction     isomorphisms:
$(\widehat{W}^{\phi})_{\phi}  \iso   (\widehat{W}_H)_{\phi}$  yields  the
correspondence of  the Theorem. \hfill $\square$

\footnotesize{

}

\footnotesize{
{\bf V.B.}: Department of Mathematics, University of Chicago,
Chicago IL
60637, USA;\\
\hphantom{x}\ab {\bf barashek@math.uchicago.edu}


{\bf S.E.}: Department of Mathematics, University of Notre Dame, Notre
Dame IN, 46656 and\\
\hphantom{x}\ab Department of Mathematics, University of Arizona,
Tucson AZ, 85721, USA;\\
\hphantom{x}\ab {\bf evens.1@nd.edu}


{\bf V.G.}: Department of Mathematics, University of Chicago,
Chicago IL
60637, USA;\\
\hphantom{x}\ab {\bf ginzburg@math.uchicago.edu}
}

\end{document}